\documentclass[a4paper,12pt]{article}
\usepackage[T1]{fontenc}
\usepackage{latexsym}
\usepackage{amssymb}
\usepackage{amsmath}
\usepackage{amsthm}
\usepackage{graphicx}
\usepackage{hyperref}
\usepackage{epsfig}
\usepackage{makeidx}
\usepackage{graphics}
\usepackage{times}
\usepackage{fancybox}
\usepackage{color}
\usepackage{colortbl}

\newtheorem{de}{Definition}[section] 
\newtheorem{coro}{Corollary} 
\newtheorem{prop}{Proposition}[section] 

\title{Differential Geometry and Mechanics\\
Applications to \\Chaotic Dynamical Systems}

\date{\today}

\author{JEAN-MARC GINOUX AND BRUNO ROSSETTO\\
\textit{PROTEE Laboratory,}
\\
\textit{I.U.T. of Toulon, University of South, }
\\
\textit{B.P. 20132, 83957, LA GARDE Cedex, France}
\\
\textcolor{blue}{ginoux@univ-tln.fr},
\textcolor{blue}{rossetto@univ-tln.fr}}

\begin{document}

\maketitle

\begin{abstract}

\textit{The aim of this article is to highlight the interest to
apply Differential} \textit{Geometry} \textit{and} \textit{Mechanics
concepts} \textit{to chaotic dynamical systems study.} \textit{Thus,
the local metric properties of curvature and torsion will directly
provide the analytical expression of the slow manifold equation of
slow-fast autonomous dynamical systems starting from kinematics
variables (velocity, acceleration and over-acceleration or jerk).}
\\
\textit{The attractivity of the slow manifold will be characterized
thanks to a criterion proposed by Henri Poincar\'{e}. Moreover, the
specific use of acceleration will make it possible on the one hand
to define slow and fast domains of the phase space and on the other
hand, to provide an analytical equation of the slow manifold towards
which all the trajectories converge. The attractive slow manifold
constitutes a part of these dynamical systems attractor. So, in
order to propose a description of the geometrical structure of
attractor, a new manifold called singular manifold will be
introduced. Various applications of this new approach to the models
of Van der Pol, cubic-Chua, Lorenz, and Volterra-Gause are
proposed.}~

\vspace{1cm}


\end{abstract}

\section{Introduction}

There are various methods to determine the \textit{slow manifold}
analytical equation of slow-fast autonomous dynamical systems
(S-FADS), or of autonomous dynamical systems considered as slow-fast
(CAS-FADS). A classical approach based on the works of Andronov
[1966] led to the famous \textit{singular approximation}. For
$\varepsilon \neq 0$, another method called: \textit{tangent linear
system approximation}, developed by Rossetto \textit{et al.} [1998],
consists in using the presence of a ``fast'' eigenvalue in the
functional jacobian matrix of a (S-FADS) or of a (CAS-FADS). Within
the framework of application of the Tihonov's theorem [1952], this
method uses the fact that in the vicinity of the \textit{slow
manifold }the eigenmode associated in the ``fast'' eigenvalue is
evanescent. Thus, the \textit{tangent linear system approximation}
method, presented in the appendix, provides the \textit{slow
manifold} analytical equation of a dynamical system according to the
``slow'' eigenvectors of the \textit{tangent linear system}, i.e.,
according to the ``slow'' eigenvalues. Nevertheless, according to
the nature of the ``slow'' eigenvalues (real or complex conjugated)
the plot of the \textit{slow manifold} analytical equation may be
difficult even impossible. Also to solve this problem it was
necessary to make the \textit{slow manifold} analytical equation
independent of the ``slow'' eigenvalues. This could be carried out
by multiplying the \textit{slow manifold} analytical equation of a
two dimensional dynamical system by a ``conjugated'' equation, that
of a three dimensional dynamical system by two ``conjugated''
equations. In each case, the \textit{slow manifold} analytical
equation independent of the ``slow'' eigenvalues of the
\textit{tangent linear system} is presented in the appendix.

The new approach proposed in this article is based on the use of
certain properties of \textit{Differential Geometry} and
\textit{Mechanics}. Thus, the metric properties of
\textit{curvature} and \textit{torsion} have provided a direct
determination of the \textit{slow manifold} analytical equation
independently of the ``slow'' eigenvalues. It has been demonstrated
that the equation thus obtained is completely identical to that
which the \textit{tangent linear system} \textit{approximation}
method provides.

The attractivity or repulsivity of the \textit{slow manifold} could
be characterized while using a criterion proposed by Henri
Poincar\'{e} [1881] in a report entitled ``Sur les courbes
d\'{e}finies par une equation diff\'{e}rentielle''.

Moreover, the specific use of the instantaneous acceleration vector
allowed a kinematic interpretation of the evolution of the
\textit{trajectory curves} in the vicinity of the \textit{slow
manifold} by defining the \textit{slow} and \textit{fast} domains of
the phase space.

At last, two new manifolds called \textit{singular} were introduced.
It was shown that the first, the \textit{singular approximation of
acceleration}, constitutes a first approximation of the \textit{slow
manifold} analytical equation of a (S-FADS) and is completely
equivalent to the equation provided by the method of the
\textit{successive approximations} developed by Rossetto [1986]
while the second, the \textit{singular manifold}, proposes an
interpretation of the geometrical structure of these dynamical
systems attractor.

\section{Dynamical System, Slow-Fast Autonomous Dynamical System (S-FADS), Considered as
Slow-Fast Autonomous Dynamical System (CAS-FADS)}

The aim of this section is to recall definitions and properties of
(S-FADS) and of (CAS-FADS)

\subsection{Dynamical system}

In the following we consider a system of differential equations
defined in a compact E included in $\mathbb{R}$:

\begin{equation}
\label{eq1} \frac{d\vec {X}}{dt} = \overrightarrow \Im \left( \vec
{X} \right)
\end{equation}

with

\[
\vec {X} = \left[ {x_1 ,x_2 ,...,x_n } \right]^t \in E \subset
\mathbb{R}^n
\]

and

\[
\overrightarrow \Im \left( \vec {X} \right) = \left[ {f_1 \left(
\vec {X} \right),f_2 \left( \vec {X} \right),...,f_n \left( \vec {X}
\right)} \right]^t \in E \subset \mathbb{R}^n
\]

The vector $\overrightarrow \Im \left( \vec {X} \right)$ defines a
velocity vector field in E whose components $f_i$ which are supposed
to be continuous and infinitely differentiable with respect to all
$x_i $ and $t$, i.e., are $C^\infty $ functions in E and with values
included in $\mathbb{R}$, satisfy the assumptions of the
Cauchy-Lipschitz theorem. For more details, see for example
Coddington {\&} Levinson [1955]. A solution of this system is an
integral curve $\vec {X}\left( t \right)$ tangent to
$\overrightarrow \Im $ whose values define the \textit{states} of
the \textit{dynamical system} described by the Eq. (\ref{eq1}).
Since none of the components $f_i $ of the velocity vector field
depends here explicitly on time, the system is said to be
\textit{autonomous}.
\\

\begin{flushleft}
\textbf{Note:}
\end{flushleft}

In certain applications, it would be supposed that the components
$f_i $ are $C^r$ functions in E and with values in $\mathbb{R}$,
with $r \geqslant n$.

\subsection{Slow-fast autonomous dynamical system (S-FADS)}

A (S-FADS) is a \textit{dynamical system} defined under the same
conditions as above but comprising a small multiplicative parameter
$\varepsilon $ in one or several components of its velocity vector
field:

\begin{equation}
\label{eq2} \frac{d\vec {X}}{dt} = \overrightarrow \Im \left( \vec
{X} \right)
\end{equation}

with

\[
\frac{d\vec {X}}{dt} = \left[ {\varepsilon \frac{dx_1
}{dt},\frac{dx_2 }{dt},...,\frac{dx_n }{dt}} \right]^t \in E \subset
\mathbb{R}^n
\]

\[
0 < \varepsilon \ll 1
\]

The functional jacobian of a (S-FADS) defined by (\ref{eq2}) has an
eigenvalue called ``fast'', i.e., with a large real part on a large
domain of the phase space. Thus, a ``fast'' eigenvalue is expressed
like a polynomial of valuation $ - 1$ in $\varepsilon $ and the
eigenmode which is associated in this ``fast'' eigenvalue is said:\\

- ``evanescent'' if it is negative,

- ``dominant'' if it is positive.\\

The other eigenvalues called ``slow'' are expressed like a
polynomial of valuation $0$ in $\varepsilon $.

\subsection{Dynamical system considered as slow-fast (CAS-FADS)}

It has been shown [Rossetto \textit{et al.}, 1998] that a dynamical
system defined under the same conditions as (\ref{eq1}) but without
small multiplicative parameter in one of the components of its
velocity vector field, and consequently without \textit{singular
approximation}, can be considered as \textit{slow}-\textit{fast} if
its functional jacobian matrix has at least a ``fast'' eigenvalue,
i.e., with a large real part on a large domain of the phase space.

\newpage

\section{New approach of the slow manifold of dynamical systems}

This approach consists in applying certain concepts of
\textit{Mechanics} and \textit{Differential Geometry} to the study
of dynamical systems (S-FADS or CAS-FADS). \textit{Mechanics} will
provide an interpretation of the behavior of the \textit{trajectory
curves}, integral of a (S-FADS) or of a (CAS-FADS), during the
various phases of their motion in terms of \textit{kinematics}
variables: velocity and acceleration. The use of
\textit{Differential Geometry}, more particularly the local metric
properties of \textit{curvature} and \textit{torsion}, will make it
possible to directly determine the analytical equation of the
\textit{slow} \textit{manifold} of (S-FADS) or of (CAS-FADS).

\subsection{Kinematics vector functions}

Since our proposed approach consists in using the \textit{Mechanics}
formalism, it is first necessary to define the \textit{kinematics}
variables needed for its development. Thus, we can associate the
integral of the system (\ref{eq1}) or (\ref{eq2}) with the
co-ordinates, i.e., with the position, of a moving point M at the
instant$ t$. This integral curve defined by the vector function
$\vec {X}\left( t \right)$ of the scalar variable $t$ represents the
\textit{trajectory curve} of the moving point M.

\subsubsection{Instantaneous velocity vector}

As the vector function $\vec {X}\left( t \right)$ of the scalar
variable $t$ represents the \textit{trajectory} of M, the total
derivative of $\vec {X}\left( t \right)$ is the vector function
$\overrightarrow V \left( t \right)$ of the scalar variable $t$
which represents the instantaneous velocity vector of the mobile M
at the instant $t$; namely:

\begin{equation}
\label{eq3} \overrightarrow V \left( t \right) = \frac{d\vec
{X}}{dt} = \overrightarrow \Im \left( \vec {X} \right)
\end{equation}

The instantaneous velocity vector $\overrightarrow V \left( t
\right)$ is supported by the tangent to the \textit{trajectory
curve}.

\newpage

\subsubsection{Instantaneous acceleration vector}

As the instantaneous vector function $\overrightarrow V \left( t
\right)$ of the scalar variable $t$ represents the velocity vector
of M, the total derivative of $\overrightarrow V \left( t \right)$
is the vector function $\overrightarrow \gamma \left( t \right)$ of
the scalar variable $t$ which represents the instantaneous
acceleration vector of the mobile M at the instant $t$; namely:

\begin{equation}
\label{eq4} \overrightarrow \gamma \left( t \right) =
\frac{d\overrightarrow V }{dt}
\end{equation}

Since the functions $f_i $ are supposed to be $C^\infty $ functions
in a compact E included in $\mathbb{R}^n$, it is possible to
calculate the total derivative of the vector field $\overrightarrow
V \left( t \right)$ defined by (\ref{eq1}) or (\ref{eq2}). By using
the derivatives of composite functions, we can write the derivative
in the sense of Fr\'{e}chet:

\begin{equation}
\label{eq5} \frac{d\overrightarrow V }{dt} = \frac{d\overrightarrow
\Im }{d\vec {X}}\frac{d\vec {X}}{dt}
\end{equation}

By noticing that $\frac{d\overrightarrow \Im }{d\vec {X}}$ is the
functional jacobian matrix $J$ of the system (\ref{eq1}) or
(\ref{eq2}), it follows from Eqs. (\ref{eq4}) and (\ref{eq5}) that
we have the following equation which plays a very important role:

\begin{equation}
\label{eq6} \overrightarrow \gamma = J\overrightarrow V
\end{equation}

\subsubsection{Tangential and normal components of the instantaneous
acceleration vector}

By making the use of the Fr\'{e}net [1847] frame, i.e., a frame
built starting from the \textit{trajectory} \textit{curve} $\vec
{X}\left( t \right)$ directed towards the motion of the mobile M.
Let's define $\overrightarrow \tau $ the unit tangent vector to the
\textit{trajectory curve} in M, $\overrightarrow \nu $ the unit
normal vector, i.e., the principal normal in M directed towards the
interior of the concavity of the curve and $\vec {\beta }$ the unit
binormal vector to the \textit{trajectory curve} in M so that the
trihedron $\left( {\overrightarrow \tau ,\overrightarrow \nu ,\vec
{\beta }} \right)$ is direct. Since the instantaneous velocity
vector
 $\overrightarrow V $ is tangent to any point M to the \textit{trajectory} \textit{curve} $\vec
{X}\left( t \right)$, we can construct a unit tangent vector as
following:

\begin{equation}
\label{eq7} \overrightarrow \tau = \frac{\overrightarrow V }{\left\|
{\overrightarrow V } \right\|}
\end{equation}

\newpage

In the same manner, we can construct a unit binormal, as:

\begin{equation}
\label{eq8} \vec {\beta } = \frac{\overrightarrow V \wedge
\overrightarrow \gamma }{\left\| {\overrightarrow V \wedge
\overrightarrow \gamma } \right\|}
\end{equation}

and a unit normal vector, as:

\begin{equation}
\label{eq9} \vec {\nu } = \vec {\beta } \wedge \vec {\tau } =
\frac{\dot {\vec {\tau }}}{\left\| \dot {\vec {\tau }} \right\|} =
\frac{\overrightarrow V ^ \bot }{\left\| {\overrightarrow V ^ \bot }
\right\|}
\end{equation}

with

\begin{equation}
\label{eq10} \left\| {\overrightarrow V } \right\| = \left\|
{\overrightarrow V ^ \bot } \right\|
\end{equation}

where the vector $\overrightarrow V ^ \bot $ represents the normal
vector to the instantaneous velocity vector $\overrightarrow V $
directed towards the interior of the concavity of the
\textit{trajectory curve} and where the dot $(\cdot)$ represents the
derivative with respect to time. Thus, we can express the tangential
and normal components of the instantaneous acceleration vector
$\overrightarrow \gamma $ as:

\begin{equation}
\label{eq11} \gamma _\tau = \frac{\overrightarrow \gamma \cdot
\overrightarrow V }{\left\| {\overrightarrow V } \right\|}
\end{equation}

\begin{equation}
\label{eq12} \gamma _\nu = \frac{\left\| {\overrightarrow \gamma
\wedge \overrightarrow V } \right\|}{\left\| {\overrightarrow V }
\right\|}
\end{equation}

By noticing that the variation of the Euclidian norm of the
instantaneous velocity vector $\overrightarrow V $ can be written:

\begin{equation}
\label{eq13} \frac{d\left\| {\overrightarrow V } \right\|}{dt} =
\frac{\overrightarrow \gamma \cdot \overrightarrow V }{\left\|
{\overrightarrow V } \right\|}
\end{equation}

And while comparing Eqs. (\ref{eq11}) and (\ref{eq12}) we deduce
that

\begin{equation}
\label{eq14} \frac{d\left\| {\overrightarrow V } \right\|}{dt} =
\gamma _\tau
\end{equation}

\newpage

Taking account of the Eq. (\ref{eq10}) and using the definitions of
the scalar and vector products, the expressions of the tangential
(\ref{eq11}) and normal (\ref{eq12}) components of the instantaneous
acceleration vector $\overrightarrow \gamma $ can be finally
written:

\begin{equation}
\label{eq15} \gamma _\tau = \frac{\overrightarrow \gamma \cdot
\overrightarrow V }{\left\| {\overrightarrow V } \right\|} =
\frac{d\left\| {\overrightarrow V } \right\|}{dt} = \left\|
{\overrightarrow \gamma } \right\|Cos\left(
{\widehat{\overrightarrow \gamma ,\overrightarrow V }} \right)
\end{equation}

\begin{equation}
\label{eq16} \gamma _\nu = \frac{\left\| {\overrightarrow \gamma
\wedge \overrightarrow V } \right\|}{\left\| {\overrightarrow V }
\right\|} = \left\| {\overrightarrow \gamma } \right\|\left|
{Sin\left( {\widehat{\overrightarrow \gamma ,\overrightarrow V }}
\right)} \right|
\end{equation}

\vspace{0.1in}

\begin{flushleft}
\textbf{Note:}
\end{flushleft}

While using the \textit{Lagrange identity}: $\left\|
{\overrightarrow \gamma \wedge \overrightarrow V } \right\|^2 +
\left( {\overrightarrow \gamma \cdot \overrightarrow V } \right)^2 =
\left\| {\overrightarrow \gamma } \right\|^2 \cdot \left\|
{\overrightarrow V } \right\|^2$, one finds easily the norm of the
instantaneous acceleration vector $\overrightarrow \gamma \left( t
\right)$.

\[
\left\| {\overrightarrow \gamma } \right\|^2 = \gamma _\tau ^2 +
\gamma _\nu ^2 = \frac{\left\| {\overrightarrow \gamma \wedge
\overrightarrow V } \right\|^2}{\left\| {\overrightarrow V }
\right\|^2} + \frac{\left( {\overrightarrow \gamma \cdot
\overrightarrow V } \right)^2}{\left\| {\overrightarrow V }
\right\|^2} = \frac{\left\| {\overrightarrow \gamma \wedge
\overrightarrow V } \right\|^2 + \left( {\overrightarrow \gamma
\cdot \overrightarrow V } \right)^2}{\left\| {\overrightarrow V }
\right\|^2} = \left\| {\overrightarrow \gamma } \right\|^2
\]

\subsection{Trajectory curve properties}

In this approach the use of \textit{Differential Geometry} will
allow a study of the metric properties of the \textit{trajectory}
\textit{curve}, i.e., \textit{curvature} and \textit{torsion} whose
definitions are recalled in this section. One will find, for
example, in Delachet [1964], Struik [1934], Kreyzig [1959] or Gray
[2006] a presentation of these concepts.

\subsubsection{Parametrization of the trajectory curve}

The \textit{trajectory curve} $\vec {X}\left( t \right)$ integral of
the dynamical system defined by (\ref{eq1}) or by (\ref{eq2}), is
described by the motion of a current point M position of which
depends on a variable parameter: the time. This curve can also be
defined by its parametric representation relative in a frame:

\[
x_1 = F_1 \left( t \right),\mbox{ }x_2 = F_2 \left( t \right),\mbox{
}...,\mbox{ }x_n = F_n \left( t \right)
\]

where the $F_{i}$ functions are continuous, $C^\infty $ functions
(or $C^{r+1}$ according to the above assumptions) in E and with
values in $\mathbb{R}$. Thus, the \textit{trajectory} \textit{curve
}$\vec {X}\left( t \right)$ integral of the dynamical system defined
by (\ref{eq1}) or by (\ref{eq2}), can be considered as a
\textit{plane curve} or as a \textit{space curve} having certain
metric properties like \textit{curvature} and \textit{torsion} which
will be defined below.

\subsubsection{Curvature of the trajectory curve}

Let's consider the \textit{trajectory curve }$\vec {X}\left( t
\right)$ having in $M$ an instantaneous velocity vector
$\overrightarrow V \left( t \right)$ and an instantaneous
acceleration vector $\overrightarrow \gamma \left( t \right)$, the
\textit{curvature}, which expresses the rate of changes of the
tangent to the \textit{trajectory curve}, is defined by:

\begin{equation}
\label{eq17} \frac{1}{\Re } = \frac{\left\| {\overrightarrow \gamma
\wedge \overrightarrow V } \right\|}{\left\| {\overrightarrow V }
\right\|^3} = \frac{\gamma _\nu }{\left\| {\overrightarrow V }
\right\|^2}
\end{equation}

where $\Re $ represents the \textit{radius of curvature}.

\vspace{0.1in}

\begin{flushleft}
\textbf{Note:}
\end{flushleft}

The location of the points where the local \textit{curvature} of the
\textit{trajectory curve} is null represents the location of the
points of analytical inflexion, i.e., the location of the points
where the normal component of the instantaneous acceleration vector
$\overrightarrow \gamma \left( t \right)$ vanishes.

\subsubsection{Torsion of the trajectory curve}

Let's consider the \textit{trajectory curve }$\vec {X}\left( t
\right)$ having in $M$ an instantaneous velocity vector
$\overrightarrow V \left( t \right)$, an instantaneous acceleration
vector $\overrightarrow \gamma \left( t \right)$, and an
instantaneous over-acceleration vector $\dot {\vec {\gamma }}$, the
\textit{torsion}, which expresses the difference between the
\textit{trajectory curve} and a \textit{plane} curve, is defined by:

\begin{equation}
\label{eq18} \frac{1}{\Im } = - \frac{\dot {\vec {\gamma }} \cdot
\left( {\overrightarrow \gamma \wedge \overrightarrow V }
\right)}{\left\| {\overrightarrow \gamma \wedge \overrightarrow V }
\right\|^2}
\end{equation}

where $\Im $ represents the \textit{radius of} \textit{torsion}.

\vspace{0.1in}

\begin{flushleft}
\textbf{Note:}
\end{flushleft}

A \textit{trajectory curve} whose local \textit{torsion} is null is
a curve whose \textit{osculating plane} is stationary. In this case,
the \textit{trajectory curve} is a plane curve.

\subsection{Application of these properties to the
determination of the slow manifold analytical equation}

In this section it will be demonstrated that the use of the local
metric properties of \textit{curvature} and \textit{torsion},
resulting from \textit{Differential Geometry}, provide the
analytical equation of the \textit{slow manifold} of a (S-FADS) or
of a (CAS-FADS) of dimension two or three. Moreover, it will be
established that the \textit{slow manifold} analytical equation thus
obtained is completely identical to that provided by the
\textit{tangent linear system} \textit{approximation }method
presented in the appendix.

\subsubsection{Slow manifold equation of a two dimensional dynamical
system}

\begin{prop}
\textit{The location of the points where the local curvature of the
trajectory curves integral of a two dimensional dynamical system
defined by (\ref{eq1}) or (\ref{eq2}) is null, provides the slow
manifold analytical equation associated in this system.}
\label{prop1}
\end{prop}

\begin{flushleft}
\textit{Analytical proof of the Proposition 3.1:}
\end{flushleft}

The vanishing condition of the \textit{curvature} provides:

\begin{equation}
\label{eq19} \frac{1}{\Re } = \frac{\left\| {\overrightarrow \gamma
\wedge \overrightarrow V } \right\|}{\left\| {\overrightarrow V }
\right\|^3} = 0\mbox{ } \Leftrightarrow \mbox{ }\overrightarrow
\gamma \wedge \overrightarrow V = \vec {0}\mbox{ } \Leftrightarrow
\mbox{ }\ddot {x}\dot {y} - \dot {x}\ddot {y} = 0
\end{equation}

By using the expression (\ref{eq6}), the co-ordinates of the
acceleration vector are written:

\[
\overrightarrow \gamma \left( {{\begin{array}{*{20}c}
 \ddot {x} \hfill \\
 \ddot {y} \hfill \\
\end{array} }} \right) = \left( {{\begin{array}{*{20}c}
 {a\dot {x} + b\dot {y}} \hfill \\
 {c\dot {x} + d\dot {y}} \hfill \\
\end{array} }} \right)
\]

The equation above is written:

\[
c\dot {x}^2 - \left( {a - d} \right)\dot {x}\dot {y} - b\dot {y}^2 =
0
\]

This equation is absolutely identical to the equation (A-27)
obtained by the \textit{tangent linear system approximation} method.

\newpage
\begin{flushleft}
\textit{Geometrical proof of the Proposition 3.1:}
\end{flushleft}

The vanishing condition of the \textit{curvature} provides:

\[
\frac{1}{\Re } = \frac{\left\| {\overrightarrow \gamma \wedge
\overrightarrow V } \right\|}{\left\| {\overrightarrow V }
\right\|^3} = 0\mbox{ } \Leftrightarrow \mbox{ }\overrightarrow
\gamma \wedge \overrightarrow V = \vec {0}
\]

The\textit{ tangent linear system approximation }makes it possible
to write that:

\[
\overrightarrow V = \alpha \overrightarrow {Y_{\lambda _1 } } +
\beta \overrightarrow {Y_{\lambda _2 } } \approx \beta
\overrightarrow {Y_{\lambda _2 } }
\]

While replacing in the expression (\ref{eq6}) we obtain:

\[
\overrightarrow \gamma = J\overrightarrow V = J\left( {\beta
\overrightarrow {Y_{\lambda _2 } } } \right) = \beta \lambda _2
\overrightarrow {Y_{\lambda _2 } } = \lambda _2 \overrightarrow V
\]

This shows that the instantaneous velocity and acceleration vectors
are collinear, which results in:

\[
\overrightarrow \gamma \wedge \overrightarrow V = \vec {0}
\]

\subsubsection{Slow manifold equation of a three dimensional
dynamical system}

\begin{prop}
\textit{The location of the points where the local torsion of the
trajectory curves integral of a three dimensional dynamical system
defined by (\ref{eq1}) or (\ref{eq2}) is null, provides the slow
manifold analytical equation associated in this
system.}\label{prop2}
\end{prop}

\begin{flushleft}
\textit{Analytical proof of the Proposition 3.2:}
\end{flushleft}

The vanishing condition of the \textit{torsion} provides:

\begin{equation}
\label{eq20} \frac{1}{\Im } = - \frac{\dot {\vec {\gamma }} \cdot
\left( {\overrightarrow \gamma \wedge \overrightarrow V }
\right)}{\left\| {\overrightarrow \gamma \wedge \overrightarrow V }
\right\|^2} = 0\mbox{ } \Leftrightarrow \mbox{ }\dot {\vec {\gamma
}} \cdot \left( {\overrightarrow \gamma \wedge \overrightarrow V }
\right) = 0
\end{equation}

The first corollary inherent in the \textit{tangent linear system}
\textit{approximation} method implies to suppose that the functional
jacobian matrix is stationary. That is to say

\[
\frac{dJ}{dt} = 0
\]

\newpage
Derivative of the expression (\ref{eq6}) provides:

\[
\dot {\vec {\gamma }} = J\frac{d\overrightarrow V }{dt} +
\frac{dJ}{dt}\overrightarrow V = J\overrightarrow \gamma +
\frac{dJ}{dt}\overrightarrow V = J^2\overrightarrow V +
\frac{dJ}{dt}\overrightarrow V \approx J^2\overrightarrow V
\]

The equation above is written:

\[
\left( {J^2\overrightarrow V } \right) \cdot \left( {\overrightarrow
\gamma \wedge \overrightarrow V } \right) = 0
\]

By developing this equation one finds term in the long term the
equation (A-34) obtained by the \textit{tangent linear system
approximation} method. The two equations are thus absolutely
identical.

\begin{flushleft}
\textit{Geometrical proof of the Proposition 3.2:}
\end{flushleft}

The\textit{ tangent linear system} \textit{approximation} makes it
possible to write that:

\[
\overrightarrow V = \alpha \overrightarrow {Y_{\lambda _1 } } +
\beta \overrightarrow {Y_{\lambda _2 } } + \delta \overrightarrow
{Y_{\lambda _3 } } \approx \beta \overrightarrow {Y_{\lambda _2 } }
+ \delta \overrightarrow {Y_{\lambda _3 } }
\]

While replacing in the expression (\ref{eq6}) we obtain:

\[
\overrightarrow \gamma = J\overrightarrow V = J\left( {\beta
\overrightarrow {Y_{\lambda _2 } } + \delta \overrightarrow
{Y_{\lambda _3 } } } \right) = \beta \lambda _2 \overrightarrow
{Y_{\lambda _2 } } + \delta \lambda _3 \overrightarrow {Y_{\lambda
_3 } }
\]

According to what precedes the over-acceleration vector is written:

\[
\dot {\vec {\gamma }} \approx J^2\overrightarrow V
\]

While replacing the velocity by its expression we have:

\[
\dot {\vec {\gamma }} \approx J^2\left( {\beta \overrightarrow
{Y_{\lambda _2 } } + \delta \overrightarrow {Y_{\lambda _3 } } }
\right) = \beta \lambda _2^2 \overrightarrow {Y_{\lambda _2 } } +
\delta \lambda _3^2 \overrightarrow {Y_{\lambda _3 } }
\]

Thus it is noticed that:

\[
\begin{array}{l}
 \overrightarrow V = \beta \overrightarrow {Y_{\lambda _2 } } + \delta
\overrightarrow {Y_{\lambda _3 } } \\
 \overrightarrow \gamma = \beta \lambda _2 \overrightarrow {Y_{\lambda _2 }
} + \delta \lambda _3 \overrightarrow {Y_{\lambda _3 } } \\
 \dot {\vec {\gamma }} = \beta \lambda _2^2 \overrightarrow {Y_{\lambda _2 }
} + \delta \lambda _3^2 \overrightarrow {Y_{\lambda _3 } } \\
 \end{array}
\]

This demonstrates that the instantaneous velocity, acceleration and
over-acceleration vectors are coplanar, which results in:

\[
\dot {\vec {\gamma }} \cdot \left( {\overrightarrow \gamma \wedge
\overrightarrow V } \right) = 0
\]

This equation represents the location of the points where
\textit{torsion} is null. The identity of the two methods is thus
established.

\vspace{0.1in}
\begin{flushleft}
\textbf{Note:}
\end{flushleft}

Main results of this study are summarized in the table (\ref{tab1})
presented below.

Abbreviations mean:

- T.L.S.A.: Tangent Linear System Approximation

- A.D.G.F.: Application of the Differential Geometry Formalism

\begin{table}[htbp]
\begin{center}
\begin{tabular}{|c|c|c|}
\hline
\multicolumn{3}{|c|}{Determination of the slow manifold analytical equation}  \\
\hline & T.L.S.A.  &  A.D.G.F. \\
\hline Dimension 2 & $\overrightarrow V \wedge \overrightarrow
{Y_{\lambda _2 } } = \vec {0}\mbox{ }$ &
$\frac{1}{\Re } = \frac{\left\| {\overrightarrow \gamma \wedge \overrightarrow V } \right\|}{\left\| {\overrightarrow V } \right\|^3} = 0\mbox{ } \Leftrightarrow \mbox{ }\overrightarrow \gamma \wedge \overrightarrow V = \vec {0}$  \\
\hline  Dimension 3 & $\overrightarrow V .\left( {\overrightarrow
{Y_{\lambda _2 } } \wedge \overrightarrow {Y_{\lambda _3 } } }
\right) = 0$ &
$\frac{1}{\Im } = - \frac{\dot {\vec {\gamma }} \cdot \left( {\overrightarrow \gamma \wedge \overrightarrow V } \right)}{\left\| {\overrightarrow \gamma \wedge \overrightarrow V } \right\|^2} = 0\mbox{ } \Leftrightarrow \mbox{ }\dot {\vec {\gamma }} \cdot \left( {\overrightarrow \gamma \wedge \overrightarrow V } \right) = 0$  \\
\hline
\end{tabular}
\end{center}
\caption{Determination of the slow manifold analytical equation}
\label{tab1}
\end{table}

In the first report entitled on the ``Courbes d\'{e}finies par une
\'{e}quation diff\'{e}rentielle'' Henri Poincar\'{e} [1881] proposed
a criterion making it possible to characterize the attractivity or
the repulsivity of a manifold. This criterion is recalled in the
next section.

\subsubsection{Attractive, repulsive manifolds}

\begin{prop}
\textit{Let }$\vec {X}\left( t \right)$ \textit{ be a trajectory
curve having in M an instantaneous velocity vector} $\overrightarrow
V \left( t \right)$ \textit{ and let ($\mathcal{V}$) be a manifold
(a curve in dimension two, a surface in dimension three) defined by
the implicit equation} $\phi = 0$ \textit{whose normal vector}
$\overrightarrow \eta = \overrightarrow \nabla \phi$ \textit{is
directed towards the outside of the concavity of this manifold.}
\\

- \textit{If the scalar product between the instantaneous velocity
vector} $\overrightarrow V \left( t \right)$ \textit{and the normal
vector} $\overrightarrow \eta = \overrightarrow \nabla \phi$
 \textit{is positive, the manifold is said attractive with respect
to this trajectory curve}

- \textit{If it is null, the trajectory curve is tangent to this
manifold.}

- \textit{If it is negative, the manifold is said repulsive.}

\textit{This scalar product which represents the total derivative
of} $\phi$ \textit{constitutes a new manifold} ($\mathcal{\dot{V}}$)
\textit{which is the envelope of the manifold ($\mathcal{V}$).}
\label{prop3}
\end{prop}

\newpage

\begin{proof}\mbox{                  }\\
\\
Let's consider a \textit{manifold} ($\mathcal{V}$) defined by the
implicit equation $\phi \left( {x,y,z} \right) = 0$.\\

The normal vector directed towards the outside of the concavity of
the \textit{curvature} of this \textit{manifold} is written:

\begin{equation}
\label{eq21} \overrightarrow \eta = \overrightarrow \nabla \phi =
\left( {{\begin{array}{*{20}c}
 {\frac{\partial \phi }{\partial x}} \hfill \\
 {\frac{\partial \phi }{\partial y}} \hfill \\
 {\frac{\partial \phi }{\partial z}} \hfill \\
\end{array} }} \right)
\end{equation}

The instantaneous velocity vector of the \textit{trajector}y
\textit{curve} is defined by (\ref{eq1}):

\[
\overrightarrow V = \left( {{\begin{array}{*{20}c}
 {\frac{dx}{dt}} \hfill \\
 {\frac{dy}{dt}} \hfill \\
 {\frac{dz}{dt}} \hfill \\
\end{array} }} \right)
\]

The scalar product between these two vectors is written:

\begin{equation}
\label{eq22} \overrightarrow V \cdot \overrightarrow \nabla \phi =
\frac{\partial \phi }{\partial x}\frac{dx}{dt} + \frac{\partial \phi
}{\partial y}\frac{dy}{dt} + \frac{\partial \phi }{\partial
z}\frac{dz}{dt}
\end{equation}

By noticing that Eq. (\ref{eq22}) represents the total derivative of
$\phi $, envelopes theory makes it possible to state that the new
\textit{manifold} ($\mathcal{\dot{V}}$) defined by this total
derivative constitutes the envelope of the \textit{manifold}
($\mathcal{V}$) defined by the equation $\phi = 0$. The
demonstration in dimension two of this Proposition results from what
precedes.
\end{proof}

\newpage

\subsubsection{Slow, fast domains}

In the \textit{Mechanics} formalism, the study of the nature of the
motion of a mobile M consists in being interested in the variation
of the Euclidian norm of its instantaneous velocity vector
$\overrightarrow V $, i.e., in the tangential component $\gamma
_\tau $ of its instantaneous acceleration vector $\overrightarrow
\gamma $. The variation of the Euclidian norm of the instantaneous
velocity vector $\overrightarrow V $ depends on the sign of the
scalar product between the instantaneous velocity vector\textbf{
}$\overrightarrow V $ and the instantaneous acceleration vector
$\overrightarrow \gamma $, i.e., the angle formed by these two
vectors. Thus if, $\overrightarrow \gamma \cdot \overrightarrow V >
0$, the variation of the Euclidian norm of the instantaneous
velocity vector\textbf{ }$\overrightarrow V $ is positive and the
Euclidian norm of the instantaneous velocity vector $\overrightarrow
V $ increases. The motion is accelerated, it is in its \textit{fast}
phase. If, $\overrightarrow \gamma \cdot \overrightarrow V = 0$, the
variation of the Euclidian norm of the instantaneous velocity vector
$\overrightarrow V $ is null and the Euclidian norm of the
instantaneous velocity vector\textbf{ }$\overrightarrow V $ is
constant. The motion is \textit{uniform}, it is in a phase of
transition between its \textit{fast} phase and its \textit{slow}
phase. Moreover, the instantaneous velocity vector $\overrightarrow
V $ is perpendicular to the instantaneous acceleration vector
$\overrightarrow \gamma $. If, $\overrightarrow \gamma \cdot
\overrightarrow V < 0$, the variation of the Euclidian norm of the
instantaneous velocity vector $\overrightarrow V $ is negative and
the Euclidian norm of the instantaneous velocity vector\textbf{
}$\overrightarrow V $ decreases. The motion is decelerated. It is in
its \textit{slow} phase.

\begin{de}

The domain of the phase space in which the tangential
component\textbf{ }$\gamma _\tau $ of the instantaneous acceleration
vector $\overrightarrow \gamma $ is negative, i.e., the domain in
which the system is decelerating is called \textit{slow} domain.

The domain of the phase space in which the tangential
component\textbf{ }$\gamma _\tau $ of the instantaneous acceleration
vector $\overrightarrow \gamma $ is positive, i.e., the domain in
which the system is accelerating is called \textit{fast} domain.

\end{de}

\begin{flushleft}
\textbf{Note:}
\end{flushleft}

On the one hand, if the (S-FADS) studied comprises only one small
multiplicative parameter $\varepsilon $ in one of the components of
its velocity vectors field, these two domains are complementary. The
location of the points belonging to the domain of the phase space
where the tangential component $\gamma _\tau $ of the instantaneous
acceleration vector $\overrightarrow \gamma $ is cancelled, delimits
the boundary between the \textit{slow} and \textit{fast} domains. On
the other hand, the slow manifold of a (S-FADS) or of a (CAS-FADS)
necessary belongs to the \textit{slow} domain.

\newpage

\section{Singular manifolds}

The use of \textit{Mechanics} made it possible to introduce a new
manifold called \textit{singular approximation of the acceleration}
which provides an approximate equation of the \textit{slow manifold}
of a (S-FADS).

\subsection{Singular approximation of the acceleration}

The\textit{ singular perturbations theory} [Andronov \textit{et
al.,} 1966] have provided the zero order approximation in
$\varepsilon $, i.e., the \textit{singular} \textit{approximation},
of the \textit{slow manifold} equation associated in a (S-FADS)
comprising a small multiplicative parameter $\varepsilon $ in one of
the components of its velocity vector field $\overrightarrow V $. In
this section, it will be demonstrated that the \textit{singular
approximation} associated in the acceleration vector field
$\overrightarrow \gamma $ constitutes the first order approximation
in $\varepsilon $ of the \textit{slow manifold} equation associated
with a (S-FADS) comprising a small multiplicative parameter
$\varepsilon $ in one of the components of its velocity vector field
$\overrightarrow V $ and consequently a small multiplicative
parameter $\varepsilon ^2$ in one of the components of its
acceleration vector field $\overrightarrow \gamma $.

\begin{prop} \textit{The manifold equation associated in
the singular approximation of the instantaneous acceleration vector
}$\overrightarrow \gamma $ \textit{of a (S-FADS) constitutes the
first order approximation in }$\varepsilon $ \textit{of the slow
manifold equation.}\label{prop4}
\end{prop}

\begin{flushleft}
\textit{Proof of the Proposition 4.1 for two-dimensional (S-FADS):}
\end{flushleft}

In dimension two, Proposition~\ref{prop1} results into a
collinearity condition (\ref{eq19}) between the instantaneous
velocity vector $\overrightarrow V $ and the instantaneous
acceleration vector $\overrightarrow \gamma $. While posing:

\[
\frac{dx}{dt} = \dot {x} \mbox{ and } \frac{dy}{dt} = \dot {y} = g
\]

The \textit{slow manifold} equation of a (S-FADS) is written:

\begin{equation}
\label{eq23} \left( {\frac{\partial g}{\partial x}} \right)\dot
{x}^2 - g\left( {\frac{1}{\varepsilon }\frac{\partial f}{\partial x}
- \frac{\partial g}{\partial y}} \right)\dot {x} - \left(
{\frac{1}{\varepsilon }\frac{\partial f}{\partial y}} \right)g^2 = 0
\end{equation}

This quadratic equation in $\dot {x}$ has the following
discriminant:

\[
\Delta = g^2\left( {\frac{1}{\varepsilon }\frac{\partial f}{\partial
x} - \frac{\partial g}{\partial y}} \right)^2 + 4\left(
{\frac{\partial g}{\partial x}} \right)\left( {\frac{1}{\varepsilon
}\frac{\partial f}{\partial y}} \right)g^2 = 0
\]

\newpage

The Taylor series of its square root up to terms of order 1 in
$\varepsilon $ is written:

\[
\sqrt \Delta \approx \frac{1}{\varepsilon }\left| {g\left(
{\frac{\partial f}{\partial x}} \right)} \right|\left\{ {1 +
\frac{\varepsilon }{\left( {\frac{\partial f}{\partial x}}
\right)^2}\left[ {2\left( {\frac{\partial g}{\partial x}}
\right)\left( {\frac{\partial f}{\partial y}} \right) - \left(
{\frac{\partial g}{\partial y}} \right)\left( {\frac{\partial
f}{\partial x}} \right)} \right] + O\left( {\varepsilon ^2} \right)}
\right\}
\]

Taking into account what precedes, the solution of the equation
(\ref{eq23}) is written:

\begin{equation}
\label{eq24} \dot {x} \approx - g\frac{\left( {\frac{\partial
f}{\partial y}} \right)}{\left( {\frac{\partial f}{\partial x}}
\right)} + O\left( \varepsilon \right)
\end{equation}

This equation represents the second-order approximation in
$\varepsilon $ of the \textit{slow manifold }equation associated
with the \textit{singular approximation}. According to equation
(\ref{eq6}), the instantaneous acceleration vector $\overrightarrow
\gamma $ is written:

\[
\overrightarrow \gamma = \left( {{\begin{array}{*{20}c}
 {\varepsilon \frac{d^2x}{dt^2}} \hfill \\
 {\frac{d^2y}{dt^2}} \hfill \\
\end{array} }} \right) = \left( {{\begin{array}{*{20}c}
 {\frac{\partial f}{\partial x}\frac{dx}{dt} + \frac{\partial f}{\partial
y}\frac{dy}{dt}} \hfill \\
 {\frac{\partial g}{\partial x}\frac{dx}{dt} + \frac{\partial g}{\partial
y}\frac{dy}{dt}} \hfill \\
\end{array} }} \right)
\]

The \textit{singular approximation of the acceleration} provides the
equation:

\[
\frac{\partial f}{\partial x}\frac{dx}{dt} + \frac{\partial
f}{\partial y}\frac{dy}{dt} = 0
\]

While posing:

\[
\frac{dx}{dt} = \dot {x} \mbox{ and } \frac{dy}{dt} = \dot {y} = g
\]

We obtain:

\begin{equation}
\label{eq25} \dot {x} = - g\frac{\left( {\frac{\partial f}{\partial
y}} \right)}{\left( {\frac{\partial f}{\partial x}} \right)}
\end{equation}

By comparing this expression with the equation (\ref{eq24}) which
constitutes the second-order approximation in $\varepsilon $ of the
\textit{slow manifold} equation, we deduces from Eq. (\ref{eq25})
that it represents the first-order approximation in $\varepsilon $
of the \textit{slow manifold} equation. It has been also
demonstrated that the \textit{singular approximation of the
acceleration} constitutes the first of the \textit{successive
approximations} developed in Rossetto [1986].

\newpage

\begin{flushleft}
\textit{Proof of the Proposition 4.1 for three-dimensional
(S-FADS):}
\end{flushleft}

In dimension three, Proposition~\ref{prop2} results into a
coplanarity condition (\ref{eq20}) between the instantaneous
velocity vector $\overrightarrow V$, the instantaneous acceleration
vector $\overrightarrow \gamma$ and the instantaneous
over-acceleration vector $\dot {\vec {\gamma }}$. The \textit{slow
manifold} equation of a (S-FADS) is written:

\begin{equation}
\label{eq26} {\frac{d^3x}{dt^3}}\left| {{\begin{array}{*{20}c}
 {\frac{dy}{dt}} \hfill & {\frac{d^2y}{dt^2}} \hfill \\
 {\frac{dz}{dt}} \hfill & {\frac{d^2z}{dt^2}} \hfill \\
\end{array} }} \right| + {\frac{d^3y}{dt^3}} \left|
{{\begin{array}{*{20}c}
 {\frac{d^2x}{dt^2}} \hfill & {\frac{dx}{dt}} \hfill \\
 {\frac{d^2z}{dt^2}} \hfill & {\frac{dz}{dt}} \hfill \\
\end{array} }} \right| + {\frac{d^3z}{dt^3}} \left|
{{\begin{array}{*{20}c}
 {\frac{dx}{dt}} \hfill & {\frac{d^2x}{dt^2}} \hfill \\
 {\frac{dy}{dt}} \hfill & {\frac{d^2y}{dt^2}} \hfill \\
\end{array} }} \right| = 0
\end{equation}

In order to simplify, let's replace the three determinants by:

\[
\Delta _1 = \left| {{\begin{array}{*{20}c}
 {\frac{dy}{dt}} \hfill & {\frac{d^2y}{dt^2}} \hfill \\
 {\frac{dz}{dt}} \hfill & {\frac{d^2z}{dt^2}} \hfill \\
\end{array} }} \right| \quad ;
\quad \Delta _2 = \left| {{\begin{array}{*{20}c}
 {\frac{d^2x}{dt^2}} \hfill & {\frac{dx}{dt}} \hfill \\
 {\frac{d^2z}{dt^2}} \hfill & {\frac{dz}{dt}} \hfill \\
\end{array} }} \right| \quad ;
\quad \Delta _3 = \left| {{\begin{array}{*{20}c}
 {\frac{dx}{dt}} \hfill & {\frac{d^2x}{dt^2}} \hfill \\
 {\frac{dy}{dt}} \hfill & {\frac{d^2y}{dt^2}} \hfill \\
\end{array} }} \right|
\]

Equation (\ref{eq26}) then will be written:

\begin{equation}
\label{eq27} \left( \dddot {x} \right)\Delta _1 + \left( \dddot {y}
\right)\Delta _2 + \left( \dddot {z} \right)\Delta _3 = 0
\end{equation}

While posing:
\[
\frac{dx}{dt} = \dot {x} \mbox{ , } \frac{dy}{dt} = \dot {y} = g
\mbox{ , } \frac{dz}{dt} = \dot {z} = h
\]

By dividing Eq. (\ref{eq27}) by $(\dddot {z})$, we have

\begin{equation}
\label{eq28} \left( {\dot {x}\ddot {y} - \ddot {x}\dot {y}} \right)
+ \frac{1}{\left( \dddot {z} \right)}\left[ {\left( \dddot {x}
\right)\Delta _1 + \left( \dddot {y} \right)\Delta _2 } \right] = 0
\end{equation}

First term of Eq. (\ref{eq28}) is written:

\begin{equation}
\label{eq29} \left( {\dot {x}\ddot {y} - \ddot {x}\dot {y}} \right)
= \left( {\frac{\partial g}{\partial x}} \right)\dot {x}^2 - g\left(
{\frac{1}{\varepsilon }\frac{\partial f}{\partial x} -
\frac{\partial g}{\partial y} - \frac{\partial g}{\partial
z}\frac{h}{g}} \right)\dot {x} - \frac{g^2}{\varepsilon }\left(
{\frac{\partial f}{\partial y} + \frac{\partial f}{\partial
z}\frac{h}{g}} \right)
\end{equation}

For homogeneity reasons let's pose:

\[
g^2G = \frac{1}{\left( \dddot {z} \right)}\left[ {\left( \dddot {x}
\right)\Delta _1 + \left( \dddot {y} \right)\Delta _2 } \right]
\]

Equation (\ref{eq28}) is written:

\begin{equation}
\label{eq30} \left( {\frac{\partial g}{\partial x}} \right)\dot
{x}^2 - g\left( {\frac{1}{\varepsilon }\frac{\partial f}{\partial x}
- \frac{\partial g}{\partial y} - \frac{\partial g}{\partial
z}\frac{h}{g}} \right)\dot {x} - \frac{g^2}{\varepsilon }\left(
{\frac{\partial f}{\partial y} + \frac{\partial f}{\partial
z}\frac{h}{g}} \right) + g^2G = 0
\end{equation}

\newpage
This quadratic equation in $\dot {x}$ has the following
discriminant:

\[
\Delta = g^2\left[ {\frac{1}{\varepsilon }\frac{\partial f}{\partial
x} - \left( {\frac{\partial g}{\partial y} + \frac{\partial
g}{\partial z}\frac{h}{g}} \right)} \right]^2 + 4g^2\left(
{\frac{\partial g}{\partial x}} \right)\left[ {\frac{1}{\varepsilon
}\left( {\frac{\partial f}{\partial y} + \frac{\partial f}{\partial
z}\frac{h}{g}} \right) - G} \right]
\]

The Taylor series of its square root up to terms of order 1 in
$\varepsilon $ is written:

\[
\sqrt \Delta \approx \frac{1}{\varepsilon }\left| {g\left(
{\frac{\partial f}{\partial x}} \right)} \right|\left\{ {1 +
\frac{\varepsilon }{\left( {\frac{\partial f}{\partial x}}
\right)^2}\left[ {2\left( {\frac{\partial g}{\partial x}}
\right)\left( {\frac{\partial f}{\partial y} + \frac{\partial
f}{\partial z}\frac{h}{g}} \right) - \left( {\frac{\partial
g}{\partial y} + \frac{\partial g}{\partial z}\frac{h}{g}}
\right)\left( {\frac{\partial f}{\partial x}} \right)} \right] +
O\left( {\varepsilon ^2} \right)} \right\}
\]

Taking into account what precedes, the solution of the equation
(\ref{eq30}) is written:

\begin{equation}
\label{eq31} \dot {x} \approx - g\frac{\left( {\frac{\partial
f}{\partial y} + \frac{\partial f}{\partial z}\frac{h}{g}}
\right)}{\left( {\frac{\partial f}{\partial x}} \right)} + O\left(
\varepsilon \right)
\end{equation}

This equation represents the second-order approximation in
$\varepsilon $ of the \textit{slow manifold }equation associated
with the \textit{singular approximation}. According to equation
(\ref{eq6}), the instantaneous acceleration vector $\overrightarrow
\gamma $ is written:

\[
\overrightarrow \gamma = \left( {{\begin{array}{*{20}c}
 {\varepsilon \frac{d^2x}{dt^2}} \hfill \\
 {\mbox{ }\frac{d^2y}{dt^2}} \hfill \\
 {\mbox{ }\frac{d^2z}{dt^2}} \hfill \\
\end{array} }} \right) = \left( {{\begin{array}{*{20}c}
 {\frac{\partial f}{\partial x}\frac{dx}{dt} + \frac{\partial f}{\partial
y}\frac{dy}{dt} + \frac{\partial f}{\partial z}\frac{dz}{dt}} \hfill \\
 {\mbox{ }\frac{\partial g}{\partial x}\frac{dx}{dt} + \frac{\partial
g}{\partial y}\frac{dy}{dt} + \frac{\partial g}{\partial
z}\frac{dz}{dt}}
\hfill \\
 {\mbox{ }\frac{\partial h}{\partial x}\frac{dx}{dt} + \frac{\partial
h}{\partial y}\frac{dy}{dt} + \frac{\partial h}{\partial
z}\frac{dz}{dt}}
\hfill \\
\end{array} }} \right)
\]

The \textit{singular approximation of the acceleration} provides the
equation:

\[
\frac{\partial f}{\partial x}\frac{dx}{dt} + \frac{\partial
f}{\partial y}\frac{dy}{dt} + \frac{\partial f}{\partial
z}\frac{dz}{dt} = 0
\]

While posing:

\[
\frac{dx}{dt} = \dot {x}, \quad \frac{dy}{dt} = \dot {y} = g \mbox{
and } \frac{dz}{dt} = \dot {z} = h
\]

We obtain:

\begin{equation}
\label{eq32} \dot {x} = - g\frac{\left( {\frac{\partial f}{\partial
y} + \frac{\partial f}{\partial z}\frac{h}{g}} \right)}{\left(
{\frac{\partial f}{\partial x}} \right)}
\end{equation}

By comparing this expression with the equation (\ref{eq31}) which
constitutes the second-order approximation in $\varepsilon $ of the
\textit{slow manifold} equation, we deduce from Eq. (\ref{eq32})
that it represents the first-order approximation in $\varepsilon $
of the \textit{slow manifold} equation.
\newpage
It has been also demonstrated that the \textit{singular
approximation of the acceleration} constitutes the first of the
\textit{successive approximations} developed in Rossetto [1986].\\

The use of the criterion proposed by H. Poincar\'{e} (Prop. 3) made
it possible to characterize the attractivity of the \textit{slow
manifold} of a (S-FADS) or of a (CAS-FADS). Moreover, the presence
in the phase space of an attractive \textit{slow} \textit{manifold},
in the vicinity of which the \textit{trajectory curves} converge,
constitutes a part of the attractor.

The\textit{ singular manifold} presented in the next section
proposes a description of the geometrical structure of the
attractor.~

\subsection{Singular manifold}

The denomination of \textit{singular manifold} comes from the fact
that this manifold plays the same role with respect to the attractor
as a singular point with respect to the \textit{trajectory curve}.

\begin{prop} \textit{The singular manifold is defined by
the intersection of slow manifold of equation }$\phi = 0$\textit{
and of an unspecified Poincar\'{e} section }$\left( \Sigma
\right)$\textit{ made in its vicinity. Thus, it represents the
location of the points satisfying:}

\begin{equation}
\label{eq33} \phi \cap \Sigma = 0
\end{equation}
\label{prop5}
\end{prop}

This \textit{manifold} of co-dimension one is a \textit{submanifold}
of the \textit{slow manifold}.\\ In dimension two, the
\textit{singular manifold} is reduced to a point.\\ In dimension
three, it is a ``line'' or more exactly a ``curve''.

The location of the points obtained by integration in a given time
of initial conditions taken on this \textit{manifold} constitutes a
\textit{submanifold} also belonging to attractor generated by the
dynamical system. The whole of these \textit{manifolds}
corresponding to various points of integration makes it possible to
reconstitute the attractor by \textit{deployment} of these
\textit{singular manifolds}.

The concept of \textit{deployment} will be illustrated in the
section 5.4.

\newpage

\section{Applications}

\subsection{Van der Pol model}

The oscillator of B. Van der Pol, [1926] is a second-order system
with non-linear frictions which can be written:

\[
\ddot {x} + \alpha (x^2 - 1)\dot {x} + x = 0
\]

The particular form of the friction which can be carried out by an
electric circuit causes a decrease of the amplitude of the great
oscillations and an increase of the small. There are various manners
of writing the previous equation like a first order system. One of
them is:

\[
\left\{ {{\begin{array}{*{20}c}
 {\dot {x} = \alpha \left( {x + y - \frac{x^3}{3}} \right)} \hfill \\
 {\dot {y} = - \frac{x}{\alpha }\mbox{ }} \hfill \\
\end{array} }} \right.
\]

When $\alpha $ becomes very large, $x$ becomes a ``fast'' variable
and $y$ a ``slow'' variable. In order to analyze the limit $\alpha
\to \infty $, we introduce a small parameter $\varepsilon = 1
\mathord{\left/ {\vphantom {1 {\alpha ^2}}} \right.
\kern-\nulldelimiterspace} {\alpha ^2}$ and a ``slow time'' $t' = t
\mathord{\left/ {\vphantom {t {\alpha = }}} \right.
\kern-\nulldelimiterspace} {\alpha = }\sqrt \varepsilon t$. Thus,
the system can be written:

\begin{equation}
\label{eq34} \overrightarrow V \left( {{\begin{array}{*{20}c}
 {\varepsilon \frac{dx}{dt}} \hfill \\
 {\mbox{ }\frac{dy}{dt}} \hfill \\
\end{array} }} \right) = \overrightarrow \Im \left( {{\begin{array}{*{20}c}
 {f\left( {x,y} \right)} \hfill \\
 {g\left( {x,y} \right)} \hfill \\
\end{array} }} \right) = \left( {{\begin{array}{*{20}c}
 {x + y - \frac{x^3}{3}} \hfill \\
 { - x} \hfill \\
\end{array} }} \right)
\end{equation}

with $\varepsilon$ a positive real parameter

\[
\varepsilon = 0.05
\]

where the functions $f$ and $g$ are infinitely differentiable with
respect to all $x_i $ and $t$, i.e., are $C^\infty $ functions in a
compact E included in $\mathbb{R}^2$ and with values in
$\mathbb{R}$. Moreover, the presence of a small multiplicative
parameter $\varepsilon $ in one of the components of its velocity
vector field $\overrightarrow V $ ensures that the system
(\ref{eq34}) is a (S-FADS). We can thus apply the method described
in Sec. 3, i.e., \textit{Differential Geometry}. The instantaneous
acceleration vector $\overrightarrow \gamma $ is written:

\begin{equation}
\label{eq35} \overrightarrow \gamma \left( {{\begin{array}{*{20}c}
 {\frac{d^2x}{dt^2}} \hfill \\
 {\frac{d^2y}{dt^2}} \hfill \\
\end{array} }} \right) = \frac{d\overrightarrow \Im }{dt}\left(
{{\begin{array}{*{20}c}
 {\frac{1}{\varepsilon }\left( {\frac{dx}{dt} + \frac{dy}{dt} -
x^2\frac{dx}{dt}} \right)} \hfill \\
 { - \frac{dx}{dt}} \hfill \\
\end{array} }} \right)
\end{equation}

Proposition~\ref{prop1} leads to:

\[
\frac{1}{\Re } = \frac{\left\| {\overrightarrow \gamma \wedge
\overrightarrow V } \right\|}{\left\| {\overrightarrow V }
\right\|^3} = 0\mbox{ } \Leftrightarrow \mbox{ }\overrightarrow
\gamma \wedge \overrightarrow V = \vec {0}\mbox{ } \Leftrightarrow
\mbox{ }\ddot {x}\dot {y} - \dot {x}\ddot {y} = 0\mbox{ }
\Leftrightarrow \mbox{ }\left| {{\begin{array}{*{20}c}
 {\frac{d^2x}{dt^2}} \hfill & {\frac{dx}{dt}} \hfill \\
 {\frac{d^2y}{dt^2}} \hfill & {\frac{dy}{dt}} \hfill \\
\end{array} }} \right| = 0
\]

We obtain the following implicit equation:

\begin{equation}
\label{eq36} \frac{1}{9\varepsilon ^2}\left[ {9y^2 + \left( {9x +
3x^3} \right)y + 6x^4 - 2x^6 + 9x^2\varepsilon } \right] = 0
\end{equation}

Since this equation is quadratic in $y$, we can solve it in order to
plot $y$ according to $x.$

\begin{equation}
\label{eq37} y_{1,2} = - \frac{x^3}{6} - \frac{x}{2}\pm
\frac{x}{2}\sqrt {x^4 - 2x^2 + 1 - 4\varepsilon }
\end{equation}

\begin{center}
\begin{figure}[htbp]
\centerline{\includegraphics{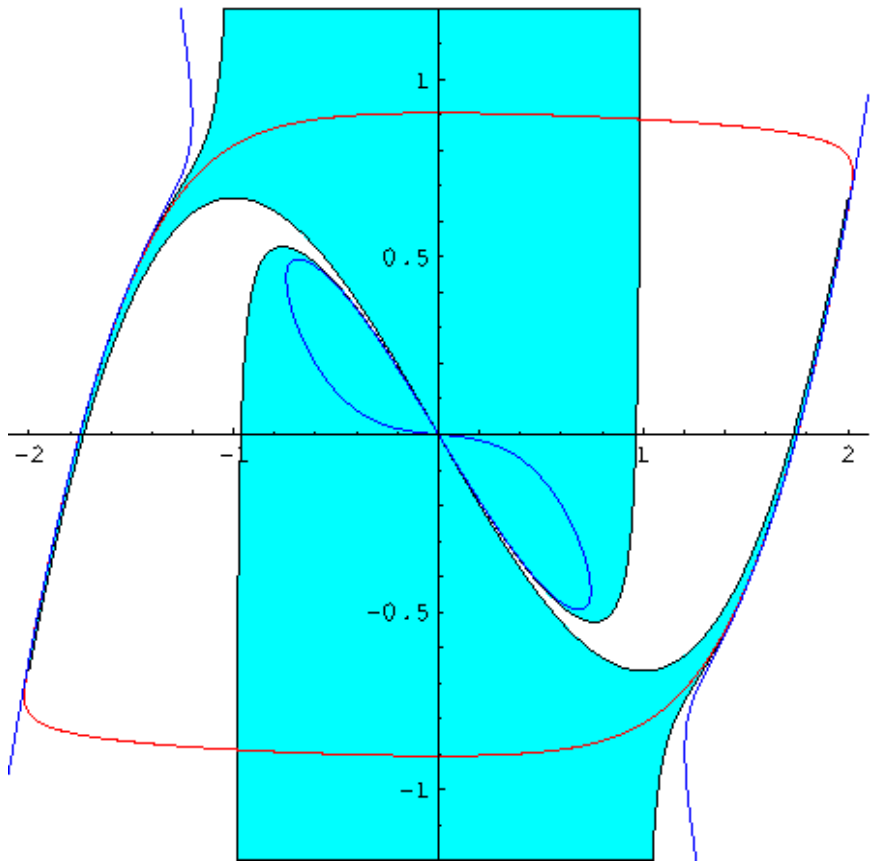}}

 \caption{Slow manifold (in blue) and stable domain of the Van der Pol system (in white).}

\label{fig1}
\end{figure}
\end{center}

In Fig.~\ref{fig1} is plotted the \textit{slow manifold} equation
(\ref{eq37}) of the Van der Pol system with $\varepsilon = 0.05$ by
using Proposition~\ref{prop1}, i.e., the collinearity condition
between the instantaneous velocity vector $\overrightarrow V $ and
the instantaneous acceleration vector $\overrightarrow \gamma $,
i.e., the location of the points where the \textit{curvature} of the
\textit{trajectory curves} is cancelled. Moreover, Definition 1
makes it possible to delimit the area of the phase plane in which,
the scalar product between the instantaneous velocity vector\textbf{
}$\overrightarrow V $ and the instantaneous acceleration
vector\textbf{ }$\overrightarrow \gamma $ is negative, i.e., where
the tangential component $\gamma _\tau $ of its instantaneous
acceleration vector $\overrightarrow \gamma $ is negative. We can
thus graphically distinguish the \textit{slow} domain of the
\textit{fast} domain (in blue), i.e., the domain of stability of the
\textit{trajectories.}

\newpage

The blue part of Fig.~\ref{fig1} corresponds to the domain where the
variation of the Euclidian norm of the instantaneous velocity vector
$\overrightarrow V $ is positive, i.e., where the tangential
component of the instantaneous acceleration vector $\overrightarrow
\gamma $, is positive. Let's notice that, as soon as the
\textit{trajectory curve}, initially outside this domain, enters
inside, it leaves the \textit{slow manifold} to reach the
\textit{fast foliation}.

\begin{center}
\begin{figure}[htbp]
\centerline{\includegraphics{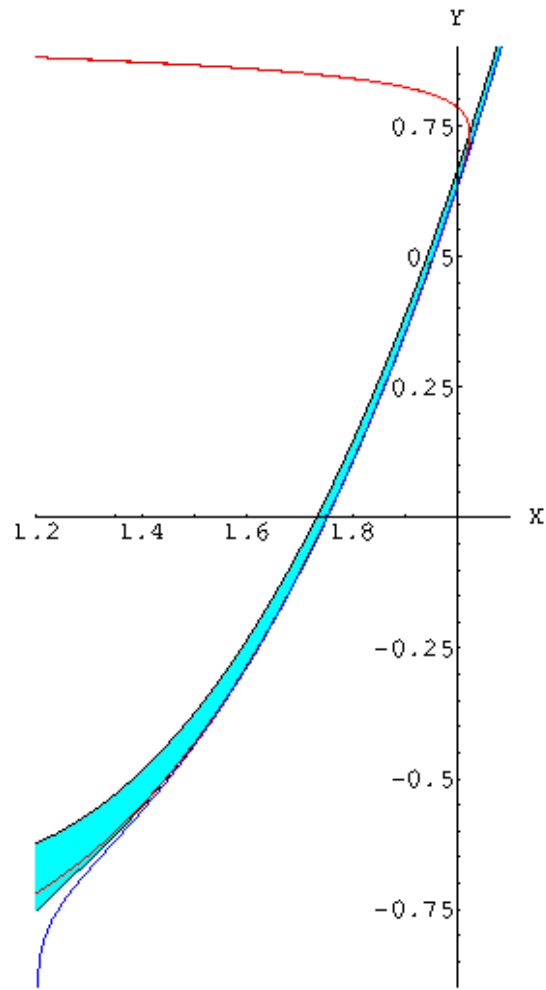}}

\caption{Zoom of the stable part of the B. Van der Pol system.}

\label{fig2}
\end{figure}
\end{center}

The \textit{slow manifold} equation provided by
Proposition~\ref{prop4} leads to the following implicit equation:

\begin{equation}
\label{eq38} \frac{1}{3\varepsilon ^2}\left[ {3x - 4x^3 + x^5 +
\left( {3 - 3x^2} \right)y - 3x\varepsilon } \right] = 0
\end{equation}

Starting from this equation we can plot $y$ according to $x$:

\begin{equation}
\label{eq39} y = \frac{x^5 - 4x^3 + 3x\left( {1 - \varepsilon }
\right)}{3\left( { - 1 + x^2} \right)}
\end{equation}

\begin{center}
\begin{figure}[htbp]
\centerline{\includegraphics{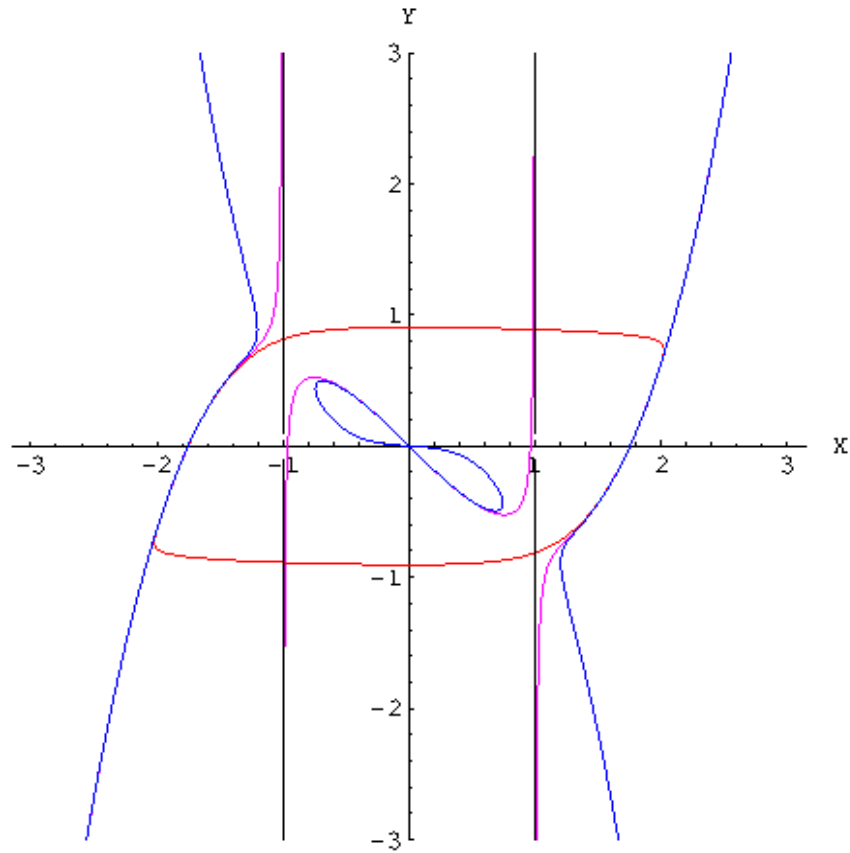}}

\caption{Singular approximation of the acceleration of the B. Van
der Pol system (in magenta).}

\label{fig3}
\end{figure}

\end{center}

In Fig.~\ref{fig3} is plotted the slow manifold equation
(\ref{eq39}) of the Van der Pol system with $\varepsilon = 0.05$ by
using Proposition~\ref{prop4}, i.e., the \textit{singular
approximation} of the instantaneous acceleration vector
$\overrightarrow \gamma$ in magenta. Blue curve represents the
\textit{slow manifold} equation (\ref{eq37}) provided by Proposition
~\ref{prop1}.

\begin{center}
\begin{figure}[htbp]
\centerline{\includegraphics{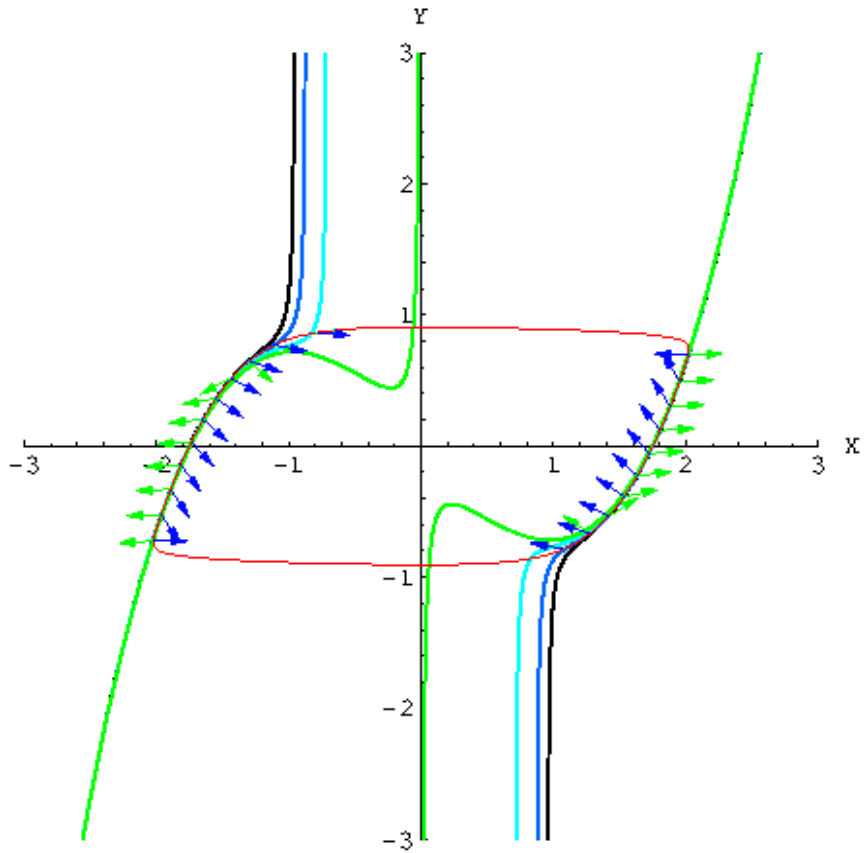}}

\caption{Isoclines of the acceleration vector of the B. Van der Pol
system.}

\label{fig4}

\end{figure}

\end{center}

In order to illustrate the principle of the method presented above,
we have plotted in the Fig.~\ref{fig4} the isoclines of acceleration
for various values: 0.5, 0.2, 0.1, 0.05.

The very large variation rate of the acceleration in the vicinity of
the \textit{slow manifold} can be noticed in Fig.~\ref{fig4} Some
isoclines of the acceleration vector which tend to the \textit{slow
manifold} defined by Proposition~\ref{prop1} are plotted.

\newpage

\subsection{Chua model}

The L.O. Chua's circuit [1986] is a relaxation oscillator with a
cubic non-linear characteristic elaborated from a circuit comprising
a harmonic oscillator of which operation is based on a field-effect
transistor, coupled to a relaxation-oscillator composed of a tunnel
diode. The modeling of the circuit uses a capacity which will
prevent from abrupt voltage drops and will make it possible to
describe the fast motion of this oscillator by the following
equations which constitute a \textit{slow-fast} system.

\begin{equation}
\label{eq40} \overrightarrow V = \left( {{\begin{array}{*{20}c}
 {\frac{dx}{dt}} \hfill \\
 {\frac{dy}{dt}} \hfill \\
 {\frac{dz}{dt}} \hfill \\
\end{array} }} \right) = \overrightarrow \Im \left( {{\begin{array}{*{20}c}
 {f\left( {x,y,z} \right)} \hfill \\
 {g\left( {x,y,z} \right)} \hfill \\
 {h\left( {x,y,z} \right)} \hfill \\
\end{array} }} \right) = \left( {{\begin{array}{*{20}c}
 {\frac{1}{\varepsilon }\left( {z - \frac{44}{3}x^3 - \frac{41}{2}x^2 - \mu
x} \right)} \hfill \\
 { - z} \hfill \\
 { - 0.7x + y + 0.24z} \hfill \\
\end{array} }} \right)
\end{equation}

with $\varepsilon$ and $\mu$ are real parameters

\[
\varepsilon = 0.05
\]

\[
\mu = 2
\]

where the functions $f$, $g$ and $h$ are infinitely differentiable
with respect to all $x_i$, and $t$, i.e., are $C^\infty$ functions
in a compact E included in $\mathbb{R}^3$ and with values in
$\mathbb{R}$. Moreover, the presence of a small multiplicative
parameter $\varepsilon$ in one of the components of its
instantaneous velocity vector $\overrightarrow V$ ensures that the
system (\ref{eq40}) is a (S-FADS). We can thus apply the method
described in Sec. 3, i.e., \textit{Differential Geometry}. In
dimension three, the \textit{slow manifold }equation is provided by
Proposition~\ref{prop2}, i.e., the vanishing condition of the
\textit{torsion}:

\[
\frac{1}{\Im } = - \frac{\dot {\vec {\gamma }} \cdot \left(
{\overrightarrow \gamma \wedge \overrightarrow V } \right)}{\left\|
{\overrightarrow \gamma \wedge \overrightarrow V } \right\|^2} =
0\mbox{ } \Leftrightarrow \mbox{ }\dot {\vec {\gamma }} \cdot \left(
{\overrightarrow \gamma \wedge \overrightarrow V } \right) = 0
\]

\begin{center}
\begin{figure}[htbp]
\centerline{\includegraphics{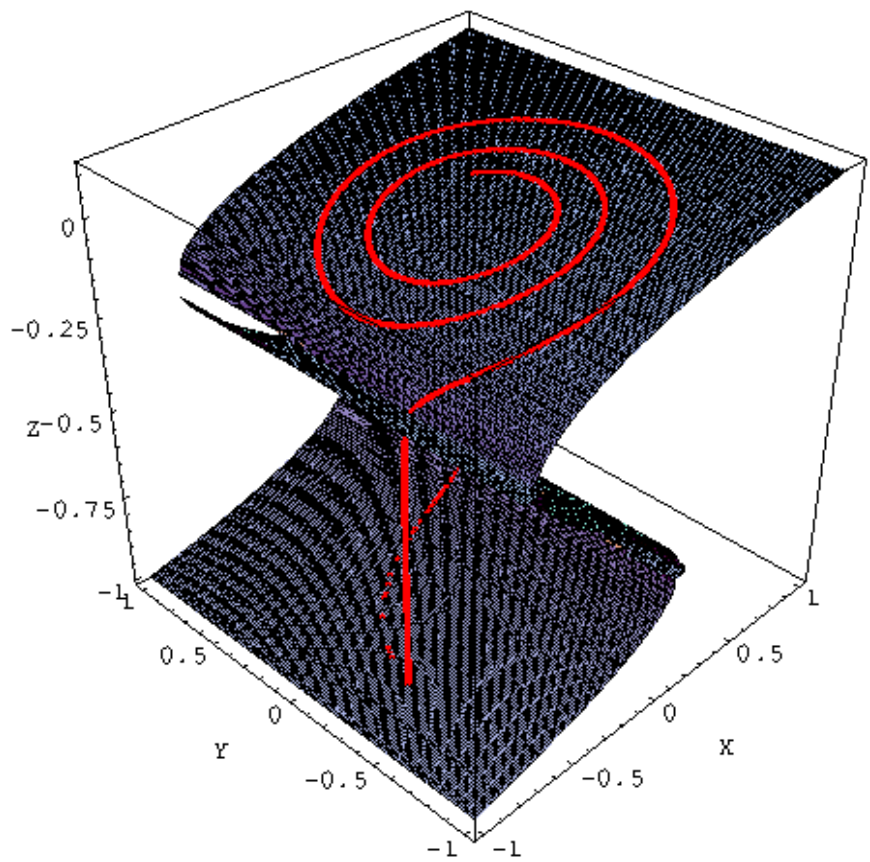}}

\caption{Slow manifold of the L.O. Chua. Model.}

\label{fig5}

\end{figure}

\end{center}

Within the framework of the \textit{tangent linear system
approximation}, Corollary 1 leads to Eq. (\ref{eq34}). While using
Mathematica$^{\circledR}$ is plotted in Fig.~\ref{fig5} the phase
portrait of the L.O. Chua model and its \textit{slow manifold. }

Without the framework of the \textit{tangent linear system
approximation}, i.e., while considering that the functional jacobian
varies with time, Proposition~\ref{prop2} provides a surface
equation which represents the location of the points where
\textit{torsion} is cancelled, i.e., the location of the points
where the \textit{osculating plane} is stationary and where the
\textit{slow manifold }is attractive.

\begin{center}
\begin{figure}[htbp]
\centerline{\includegraphics{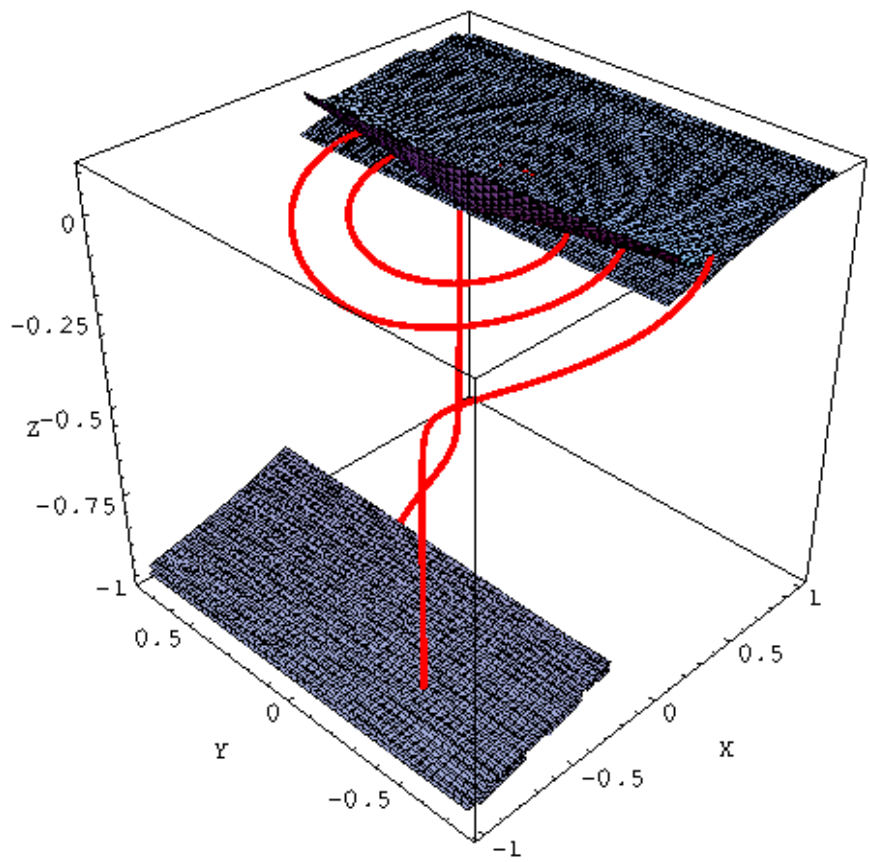}}

\caption{Attractive part of the \textit{slow manifold} of the L.O.
Chua model.}

\label{fig6}
\end{figure}

\end{center}

We deduce, according to Proposition~\ref{prop3}, that the location
of the points where the \textit{torsion} is negative corresponds to
the attractive parts of the slow manifold. Thus, the attractive part
of the \textit{slow manifold} of the LO. Chua model is plotted in
Fig.~\ref{fig6}.

\textit{Slow manifold} equation provided by Proposition~\ref{prop4}
leads to the following implicit equation:

\[
\begin{array}{l}
 \frac{1}{6\varepsilon ^2}(5043x^3 + 9020x^4 + 3872x^5 - 246xz - 264x^2z -
4.2x\varepsilon \\
 + 6y\varepsilon \mbox{ } + 1.44z\varepsilon + 369x^2\mu + 352x^3\mu - 6z\mu
+ 6x\mu ^2) = 0 \\
 \end{array}
\]

\begin{center}
\begin{figure}[htbp]
\centerline{\includegraphics{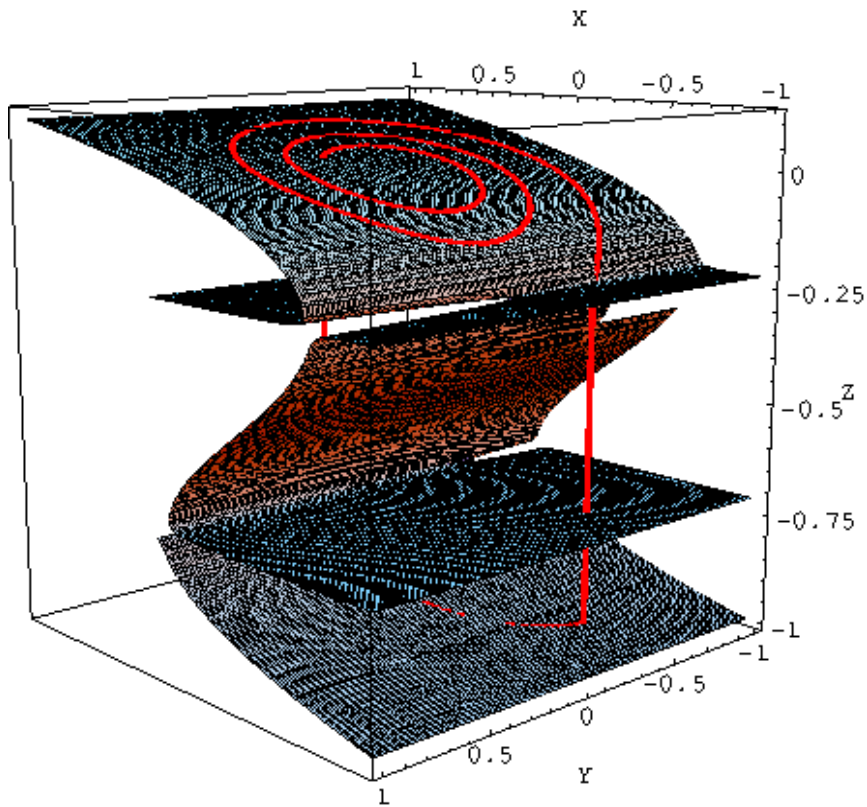}}

\caption{Singular approximation of the acceleration of the L.O. Chua
model.}

\label{fig7}
\end{figure}

\end{center}

The surface plotted in Fig.~\ref{fig7} constitutes a quite good
approximation of the \textit{slow manifold} of this model.

\section{Lorenz model}

The purpose of the model established by Edward Lorenz [1963] was in
the beginning to analyze the impredictible behavior of weather.
After having developed non-linear partial derivative equations
starting from the thermal equation and Navier-Stokes equations,
Lorenz truncated them to retain only three modes. The most
widespread form of the Lorenz model is as follows:

\begin{equation}
\label{eq41} \overrightarrow V = \left( {{\begin{array}{*{20}c}
 {\frac{dx}{dt}} \hfill \\
 {\frac{dy}{dt}} \hfill \\
 {\frac{dz}{dt}} \hfill \\
\end{array} }} \right) = \overrightarrow \Im \left( {{\begin{array}{*{20}c}
 {f\left( {x,y,z} \right)} \hfill \\
 {g\left( {x,y,z} \right)} \hfill \\
 {h\left( {x,y,z} \right)} \hfill \\
\end{array} }} \right) = \left( {{\begin{array}{*{20}c}
 {\sigma \left( {y - x} \right)} \hfill \\
 { - xz + rx - y} \hfill \\
 {xy - \beta z} \hfill \\
\end{array} }} \right)
\end{equation}

with $\sigma$, r, and $\beta $ are real parameters: $\sigma = 10$,
$\beta = \frac{8}{3}$, r = 28

where the functions $f$, $g $ and$ h$ are infinitely differentiable
with respect to all $x_i$, and $t$, i.e., are $C^\infty $ functions
in a compact E included in $\mathbb{R}^3$ and with values in
$\mathbb{R}$. Although this model has no \textit{singular
approximation}, it can be considered as a (S-FADS), according to
Sec. 1.3, because it has been numerically checked [Rossetto
\textit{et al.}, 1998] that its functional jacobian matrix possesses
at least a large and negative real eigenvalue in a large domain of
the phase space. Thus, we can apply the method described in Sec. 3,
i.e., \textit{Differential Geometry}. In dimension three, the
\textit{slow manifold }equation is provided by
Proposition~\ref{prop2}, i.e., the vanishing condition of the
\textit{torsion}:

\[
\frac{1}{\Im } = - \frac{\dot {\vec {\gamma }} \cdot \left(
{\overrightarrow \gamma \wedge \overrightarrow V } \right)}{\left\|
{\overrightarrow \gamma \wedge \overrightarrow V } \right\|^2} =
0\mbox{ } \Leftrightarrow \mbox{ }\dot {\vec {\gamma }} \cdot \left(
{\overrightarrow \gamma \wedge \overrightarrow V } \right) = 0
\]

\begin{center}
\begin{figure}[htbp]
\centerline{\includegraphics{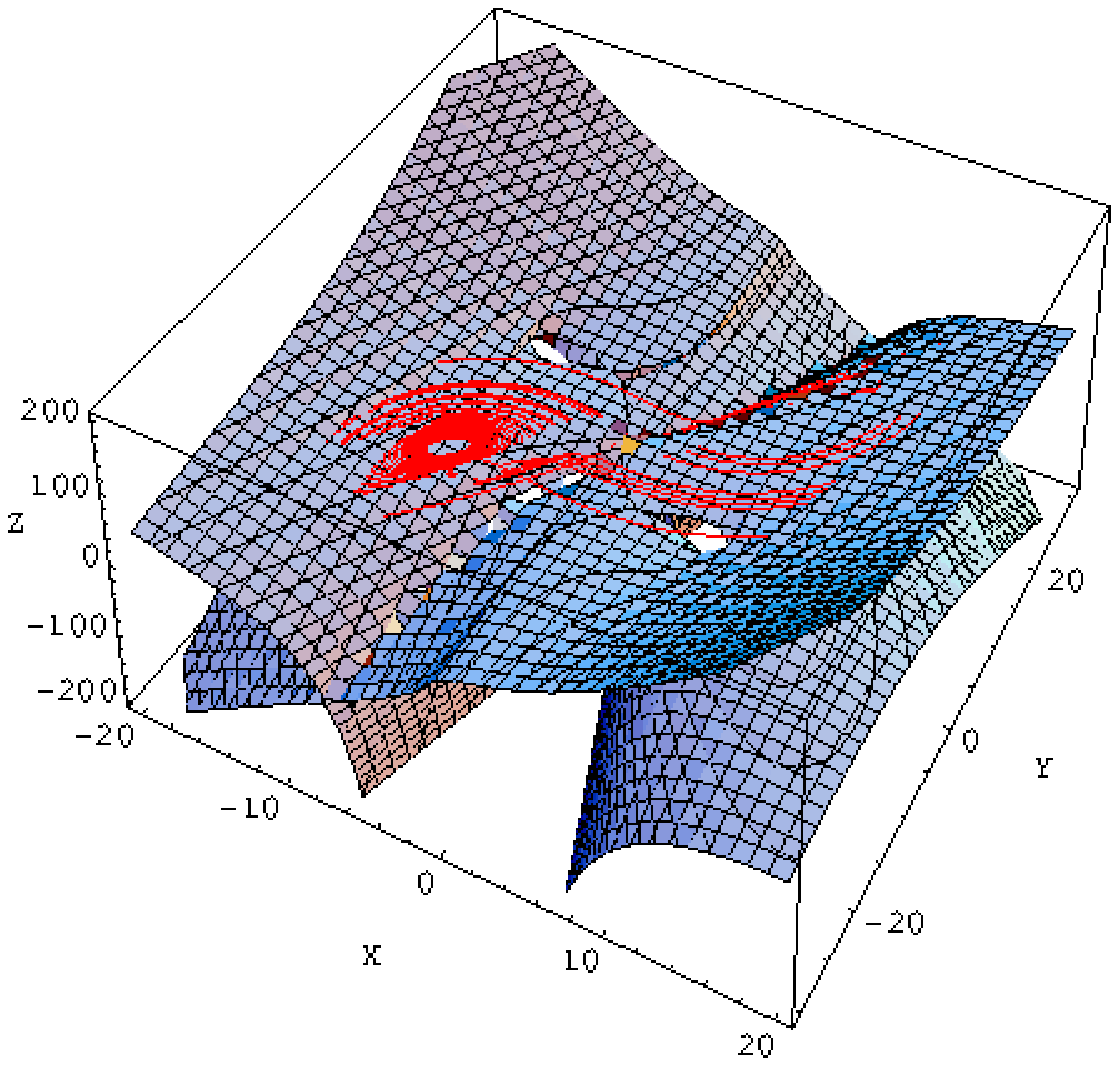}}

\caption{Slow manifold of the E.N. Lorenz model.}

\label{fig8}
\end{figure}

\end{center}

Within the framework of the \textit{tangent linear system
approximation}, Corollary 1 leads to Eq. (\ref{eq34}). While using
Mathematica$^{\circledR}$is plotted in Fig.~\ref{fig8} the phase
portrait of the E.N. Lorenz model and its \textit{slow manifold. }

\newpage
Without the framework of the \textit{tangent linear system
approximation}, i.e., while considering that the functional jacobian
varies with time, Proposition~\ref{prop2} provides a surface
equation which represents the location of the points where
\textit{torsion} is cancelled, i.e., the location of the points
where the \textit{osculating plane} is stationary and where the
\textit{slow manifold }is attractive.

\begin{center}
\begin{figure}[htbp]
\centerline{\includegraphics{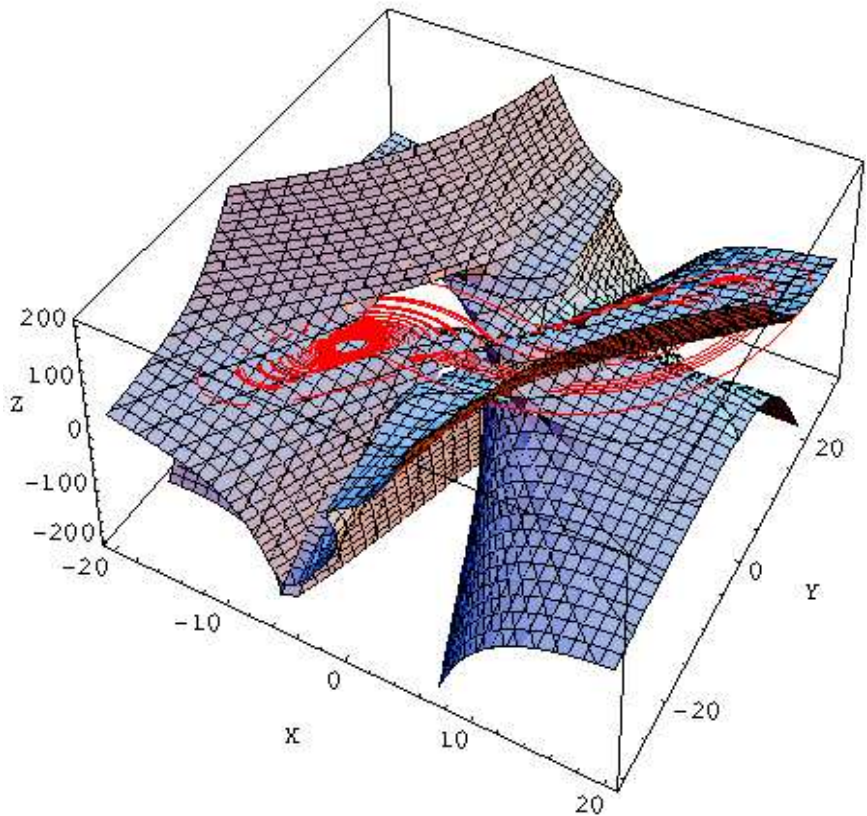}}

\caption{Attractive part of the \textit{slow manifold} of the E.N.
Lorenz model.}

\label{fig9}
\end{figure}

\end{center}

We deduce, according to Proposition~\ref{prop2}, that the location
of the points where the \textit{torsion} is negative corresponds to
the attractive parts of the \textit{slow manifold}. Thus, the
attractive part of the slow manifold of the E.N. Lorenz model is
plotted Fig.~\ref{fig9}.

\newpage

\section{Volterra-Gause model}

Let's consider the model elaborated by Ginoux \textit{et al.} [2005]
for three species interacting in a predator-prey mode.

\begin{equation}
\label{eq42} \overrightarrow V = \left( {{\begin{array}{*{20}c}
 {\frac{dx}{dt}} \hfill \\
 {\frac{dy}{dt}} \hfill \\
 {\frac{dz}{dt}} \hfill \\
\end{array} }} \right) = \overrightarrow \Im \left( {{\begin{array}{*{20}c}
 {f\left( {x,y,z} \right)} \hfill \\
 {g\left( {x,y,z} \right)} \hfill \\
 {h\left( {x,y,z} \right)} \hfill \\
\end{array} }} \right) = \left( {{\begin{array}{*{20}c}
 {\frac{1}{\xi }\left( {x\left( {1 - x} \right) - x^{\frac{1}{2}}y} \right)}
\hfill \\
 { - \delta _1 y + x^{\frac{1}{2}}y - y^{\frac{1}{2}}z} \hfill \\
 {\varepsilon z\left( {y^{\frac{1}{2}} - \delta _2 } \right)} \hfill \\
\end{array} }} \right)
\end{equation}

with $\xi$, $\varepsilon $, $\delta _1$ and $\delta _2$ are real
parameters: $\xi = 0.866$, $\varepsilon = 1.428$, $\delta _1 =
0.577$, $\delta _2 = 0.376$.

And where the functions $f, g $ and $h$ are infinitely
differentiable with respect to all $x_i $, and $t$, i.e., are
$C^\infty $ functions in a compact E included in $\mathbb{R}^3$ and
with values in $\mathbb{R}$.

\begin{center}
\begin{figure}[htbp]
\centerline{\includegraphics{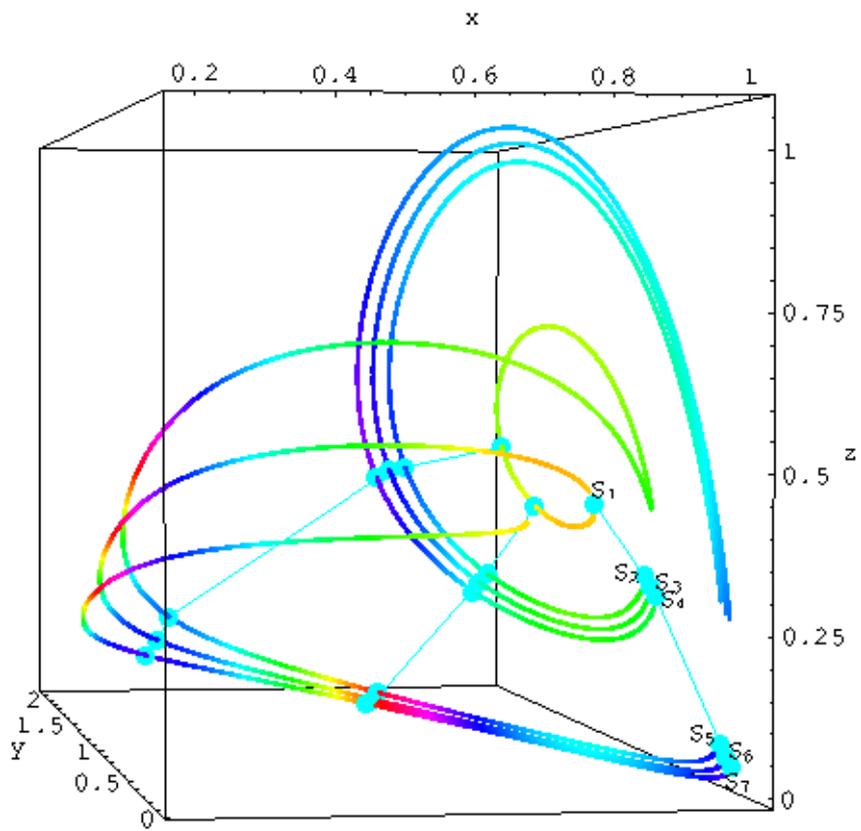}}

\caption{Deployment of the \textit{singular manifold} of the
Volterra-Gause model.}

\label{fig10}
\end{figure}

\end{center}

This model consisted of a prey, a predator and top-predator has been
named Volterra-Gause because it combines the original model of V.
Volterra [1926] incorporating a logisitic limitation of P.F.
Verhulst [1838] type on the growth of the prey and a limitation of
G.F. Gause [1935] type on the intensity of the predation of the
predator on the prey and of top-predator on the predator. The
equations (\ref{eq42}) are dimensionless, remarks and details about
the changes of variables and the parameters have been extensively
made in Ginoux \textit{et al.} [2005]. Moreover, the presence of a
small multiplicative parameter $\xi $ in one of the components of
its instantaneous velocity vector $\overrightarrow V $ ensures that
the system (\ref{eq42}) is a (S-FADS). So, the method described in
Sec. 3, i.e., \textit{Differential Geometry} would have provided the
\textit{slow manifold} equation thanks to Proposition~\ref{prop2}.
But, this model exhibits a chaotic attractor in the snail shell
shape and the use of the algorithm developed by Wolf \textit{et al.
}[1985] made it possible to compute what can be regarded as its
Lyapunov exponents: $\left( { + 0.035,0.000, - 0.628} \right)$.
Then, the Kaplan-Yorke [1983] conjecture provided the following
Lyapunov dimension: 2.06. So, the fractal dimension of this chaotic
attractor is close to that of a surface. The\textit{ singular
manifold} makes it possible to account for the evolution of the
\textit{trajectory curves} on the surface generated by this
attractor. Indeed, the location of the points of intersection of the
\textit{slow manifold} with a Poincar\'{e} section carried out in
its vicinity constitutes a ``line'' or more exactly a ``curve''.
Then by using numerical integration this ``curve'' (resp. ``line'')
is deployed through the phase space and its deployment reconstitutes
the attractor shape. The result is plotted in Fig.~\ref{fig10}.

\section{Discussion}

Considering the \textit{trajectory curves} integral of dynamical
systems as \textit{plane} or \textit{space} curves evolving in the
phase space, it has been demonstrated in this work that the local
metric properties of \textit{curvature} and \textit{torsion} of
these \textit{trajectory curves} make it possible to directly
provide the analytical equation of the \textit{slow manifold} of
dynamical systems (S-FADS or CAS-FADS) according to kinematics
variables. The \textit{slow manifold} analytical equation is
thus given by the\\

- vanishing condition of the \textit{curvature} in dimension two,

- vanishing condition of the \textit{torsion} in dimension three.\\

Thus, the use of \textit{Differential Geometry} concepts made it
possible to make the analytical equation of the \textit{slow
manifold} completely independent of the ``slow'' eigenvectors of the
functional jacobian of the \textit{tangent linear system} and it was
demonstrated that the equation thus obtained is completely identical
to that which provides the \textit{tangent linear system}
\textit{approximation} method [Rossetto \textit{et al.}, 1998]
presented below in the appendix. This made it possible to
characterize its attractivity while using a criterion proposed by
Henri Poincar\'{e} [1881] in his report entitled ``Sur les courbes
d\'{e}finies par une equation diff\'{e}rentielle''.

Moreover, the specific use of the instantaneous acceleration vector,
inherent in \textit{Mechanics}, allowed on the one hand a kinematics
interpretation of the evolution of the \textit{trajectory curves} in
the vicinity of the \textit{slow manifold} by defining the
\textit{slow} and \textit{fast} domains of the phase space and on
the other hand, to approach the analytical equation of the
\textit{slow manifold} thanks to the \textit{singular approximation
of acceleration}. The equation thus obtained is completely identical
to that which provides the \textit{successive approximations }method
[Rossetto, 1986]. Thus, it has been established that the presence in
the phase space of an attractive \textit{slow manifold}, in the
vicinity of which the \textit{trajectory curves} converge, defines
part of the attractor. So, in order to propose a qualitative
description of the geometrical structure of attractor a new manifold
called \textit{singular} has been introduced.

Various applications to the models of Van der Pol, cubic-Chua,
Lorenz and Volterra-Gause made it possible to illustrate the
practical interest of this new approach for the dynamical systems
(S-FADS or CAS-FADS) study.

\section*{Acknowledgements}

Authors would like to thank Dr. C. G\'{e}rini for his
bibliographical and mathematical remarks and comments.

\section*{R\'{e}f\'{e}rences}

Andronov, A. A., Khaikin, S. E. {\&} Vitt, A. A. [1966]
\textit{Theory of oscillators}, Pergamon Press, Oxford.
\\

Chlouverakis, K.E. {\&} Sprott, J. C. [2004]
``\href{http://plasma.physics.wisc.edu/mst/html/journalpublications/pdf/journal/Chlouverakis123.pdf}{A
Comparison of Correlation and Lyapunov Dimensions},''
\textit{Physica D}, Vol. 200, 156-164.
\\

Chua, L. O., Komuro, M. {\&} Matsumoto, T. [1986] ``The double
scroll family,'' \textit{IEEE Trans. on Circuits and Systems},
\textbf{33}, Vol.11, 1072-1097.\\

Coddington, E.A. {\&} Levinson$.$, N. [1955] \textit{Theory of
Ordinary Differential Equations}, Mac Graw Hill, New York.\\

Delachet, A. [1964] \textit{La G\'{e}om\'{e}trie
Diff\'{e}rentielle},
``\href{http://www.bibliopoche.com/collection/Que_sais-je/424.html}{Que
sais-je}'', n\r{ }1104,
\href{http://www.bibliopoche.com/editeur/PUF/104.html}{PUF}, Paris.\\

Frenet, F. [1852] ``Sur les courbes \`{a} double
courbure,''~Th\`{e}se Toulouse, 1847. R\'{e}sum\'{e} dans \textit{J.
de Math.}, 17.\\

Gause, G.F. [1935] \textit{The struggle for existence}, Williams and
Wilkins, Baltimore.\\

Ginoux, J.M., Rossetto, B. {\&} Jamet, J.L. [2005] ``Chaos in a
three-dimensional Volterra-Gause model of predator-prey
type,''\textit{ Int. J. Bifurcation and Chaos}, \textbf{5}, Vol. 15,
1689-1708.\\

Gray, A., Salamon, S. {\&} Abbena E. [2006] \textit{Modern
Differential Geometry of Curves And Surfaces With Mathematica},
Chapman {\&} Hall/CRC.\\

Kaplan, J. {\&} Yorke, J. A. [1979] ``Chaotic behavior of
multidimensional difference equations, in functional differential
equations and approximation of fixed points,'' Lecture Notes in
Mathematics, 730, 204-227.\\

Kreyszig, E. [1991]
\href{http://www.amazon.com/exec/obidos/ASIN/0486667219/weisstein-20}{Differential
Geometry,} Dover, New York.\\

Lorenz, E. N. [1963] ``Deterministic non-periodic flows,'' J. Atmos.
Sc, 20, 130-141.\\

Poincar\'{e}, H. [1881] ``Sur les courbes d\'{e}finies par une
\'{e}quation diff\'{e}rentielle,'' J. Math. Pures et Appl.,
S\'{e}rie III, 7, 375-422.\\

Poincar\'{e}, H. [1882] ``Sur les courbes d\'{e}finies par une
\'{e}quation diff\'{e}rentielle,'' J. de Math \newline Pures Appl.,
S\'{e}rie III, 8, 251-296.\\

Poincar\'{e}, H. [1885] ``Sur les courbes d\'{e}finies par une
\'{e}quation diff\'{e}rentielle,'' J. Math. Pures et Appl.,
S\'{e}rie IV, 1, 167-244.\\

Poincar\'{e}, H. [1886] ``Sur les courbes d\'{e}finies par une
\'{e}quation diff\'{e}rentielle,'' J. Math. Pures et Appl.,
S\'{e}rie IV, 2, 151-217.\\

Rossetto, B. [1986] ``Singular approximation of chaotic slow-fast
dynamical systems,'' Lecture Notes in Physics, Vol. 278, 12-14.\\

Rossetto, B., Lenzini, T., Ramdani, S. {\&} Suchey, G. [1998]
``Slow-fast autonomous dynamical systems,'' Int. J. Bifurcation and
Chaos, 8, Vol. 11, 2135-2145.\\

Struik, D.J. [1988] \textit{Lecture on Classical Differential
Geometry}, Dover, New York.\\

Tihonov, A. N. [1952] ``Systems of differential equations containing
small parameters in the derivatives,'' Mat. Sbornik N.S. 31,
575-586.\\

Van der Pol, B. [1926] ``On 'Relaxation-Oscillations','' Phil. Mag.,
7, Vol. 2, 978-992.\\

Verhulst, P.F. [1838] ``Notice sur la loi que suit la population
dans son accroissement,'' Corresp. Math. Phys., X, 113-121.\\

Volterra, V. [1926] ``Variazioni e fluttuazioni del numero
d'individui in specie animali conviventi,'' Mem. Acad. Lincei III,
6, 31-113.\\

Volterra, V. [1931] \textit{Le\c{c}ons sur la Th\'{e}orie
Math\'{e}matique de la Lutte pour la Vie}, Gauthier-Villars, Paris.\\

Wolf, A., Swift, J. B., Swinney, H. L. {\&} Vastano, J. A. [1985]
``Determining Lyapunov Exponents from a Time Series,'' Physica D,
Vol. 16, 285-317.

\section*{Appendix: formalization of the \textit{tangent
linear system approximation} method}

The aim of this appendix is to demonstrate that the approach
developed in this work generalizes the \textit{tangent linear
system} \textit{approximation }method [Rossetto \textit{et al.},
1998]. After having pointed out the necessary assumptions to the
application of this method and the corollaries which result from
this, two conditions (of collinearity / coplanarity and
orthogonality) providing the analytical equation of the \textit{slow
manifold} of a dynamical system defined by (\ref{eq1}) or
(\ref{eq2}) will be presented in a formal way. Equivalence between
these two conditions will be then established. Lastly, while using
the sum and the product (also the square of the sum and the product)
of the eigenvalues of the functional jacobian of the \textit{tangent
linear system}, the equation of the \textit{slow manifold} generated
by these two conditions will be made independent of these
eigenvalues and will be expressed according to the elements of the
functional jacobian matrix of the \textit{tangent linear system}. It
will be thus demonstrated that this analytical equation of the
\textit{slow manifold} is completely identical to that provided by
the Proposition~\ref{prop1} and ~\ref{prop2} developed in this
article.~

\subsection*{Assumptions}

The application of the \textit{tangent linear system}
\textit{approximation }method requires that the dynamical system
defined by (\ref{eq1}) or (\ref{eq2}) satisfies the following
assumptions:\\

\textbf{(H}$_{1}$\textbf{)} The components $f_i $, of the velocity
vector field $\overrightarrow \Im \left( \vec {X} \right)$ defined
in E are continuous, $C^\infty $ functions in E and with values
included in $\mathbb{R}$.\\

\textbf{(H}$_{2}$\textbf{)} The dynamical system defined by
(\ref{eq1}) or (\ref{eq2}) satisfies the \textit{nonlinear part
condition} [Rossetto \textit{et al.}, 1998], i.e., that the
influence of the nonlinear part of the Taylor series of the velocity
vector field $\overrightarrow \Im \left( \vec {X} \right)$ of this
system is overshadowed by the fast dynamics of the linear part.

\[
\overrightarrow \Im \left( \vec {X} \right) = \overrightarrow \Im
\left( {\vec {X}_0 } \right) + \left( {\vec {X} - \vec {X}_0 }
\right)\left. {\frac{d\overrightarrow \Im \left( \vec {X}
\right)}{d\vec {X}}} \right|_{\vec {X}_0 } + O\left( {\left( {\vec
{X} - \vec {X}_0 } \right)^2} \right) \qquad (A-1)
\]

\newpage

\subsection*{Corollaries}

To the dynamical system defined by (\ref{eq1}) or (\ref{eq2}) is
associated a \textit{tangent linear system }defined as follows:

\[
\frac{d\delta \vec {X}}{dt} = J\left( {\vec {X}_0 } \right)\delta
\vec {X} \qquad (A-2)
\]

where

\[
\delta \vec {X} = \vec {X} - \vec {X}_0 , \quad \vec {X}_0 = \vec
{X}\left( {t_0 } \right)\mbox{ and } \left. {\frac{d\overrightarrow
\Im \left( \vec {X} \right)}{d\vec {X}}} \right|_{\vec {X}_0 } =
J\left( {\vec {X}_0 } \right)
\]

\begin{coro}

The \textit{nonlinear part condition} implies the stability of the
\textit{slow manifold}. So, according to Proposition~\ref{prop3} it
results that the velocity varies slowly on the \textit{slow
manifold}. This involves that the functional jacobian $J\left( {\vec
{X}_0 } \right)$ varies slowly with time, i.e.,

\[
\frac{dJ}{dt}\left( {\vec {X}_0 } \right) = 0 \qquad (A-3)
\]

The solution of the \textit{tangent linear} \textit{system} (A-2) is
written:

\[
\delta \vec {X} = e^{J\left( {\vec {X}_0 } \right)\left( {t - t_0 }
\right)}\delta \vec {X}\left( {t_0 } \right) \qquad (A-4)
\]

So,

\[
\delta \vec {X} = \sum\limits_{i = 1}^n {a_i } \overrightarrow
{Y_{\lambda _i } } \qquad (A-5)
\]

where $n$ is the dimension of the eigenspace, $a_{i}$ represents
coefficients depending explicitly on the co-ordinates of space and
implicitly on time and $\overrightarrow {Y_{\lambda _i } } $ the
eigenvectors associated in the functional jacobian of the
\textit{tangent linear system}.

\end{coro}

\begin{coro}

In the vicinity of the \textit{slow manifold }the velocity of the
dynamical system defined by (\ref{eq1}) or (\ref{eq2}) and that of
the \textit{tangent linear system} (\ref{eq4}) merge.

\[
\frac{d\delta \vec {X}}{dt} = \overrightarrow V _T \approx
\overrightarrow V \qquad (A-6)
\]

where $\overrightarrow V _T $ represents the velocity vector
associated in the \textit{tangent linear system}.

\end{coro}

The \textit{tangent linear system approximation} method consists in
spreading the velocity vector field $\overrightarrow V $ on the
eigenbasis associated in the eigenvalues of the functional jacobian
of the \textit{tangent linear system}.

Indeed, by taking account of (A-2) and (A-5) we have according to
(A-6):

\[
\frac{d\delta \vec {X}}{dt} = J\left( {\vec {X}_0 } \right)\delta
\vec {X} = J\left( {\vec {X}_0 } \right)\sum\limits_{i = 1}^n {a_i }
\overrightarrow {Y_{\lambda _i } } = \sum\limits_{i = 1}^n {a_i }
J\left( {\vec {X}_0 } \right)\overrightarrow {Y_{\lambda _i } } =
\sum\limits_{i = 1}^n {a_i } \lambda _i \overrightarrow {Y_{\lambda
_i } } \qquad (A-7)
\]

Thus, Corollary 2 provides:

\[
\frac{d\delta \vec {X}}{dt} = \overrightarrow V _T \approx
\overrightarrow V = \sum\limits_{i = 1}^n {a_i } \lambda _i
\overrightarrow {Y_{\lambda _i } } \qquad (A-8)
\]

The equation (A-8) constitutes in dimension two (resp. dimension
three) a condition called \textit{collinearity} (resp.
\textit{coplanarity}) condition which provides the analytical
equation of the \textit{slow manifold} of a dynamical system defined
by (\ref{eq1}) or (\ref{eq2}).

An alternative proposed by Rossetto \textit{et al.} [1998] uses the
``fast'' eigenvector on the left associated in the ``fast''
eigenvalue of the transposed functional jacobian of the
\textit{tangent linear system}.

In this case the velocity vector field $\overrightarrow V $ is then
orthogonal with the ``fast'' eigenvector on the left. This
constitutes a condition called \textit{orthogonality} condition
which provides the analytical equation of the \textit{slow manifold}
of a dynamical system defined by (\ref{eq1}) or (\ref{eq2}).

These two conditions will be the subject of a detailed presentation
in the following sections. Thereafter it will be supposed that the
assumptions (H1) - (H2) are always checked.

\newpage

\subsection*{Collinearity / coplanarity condition}

\subsubsection*{Slow manifold equation of a two dimensional
dynamical system}

Let's consider a dynamical system defined under the same conditions
as (\ref{eq1}) or (\ref{eq2}). The eigenvectors associated in the
eigenvalues of the functional jacobian of the \textit{tangent linear
system} are written:

\[
\overrightarrow {Y_{\lambda _i } } \left( {{\begin{array}{*{20}c}
 {\lambda _i - \frac{\partial g}{\partial y}} \hfill \\
 {\frac{\partial g}{\partial x}} \hfill \\
\end{array} }} \right) \qquad (A-9)
\]

with

\begin{center}
i =1, 2
\end{center}

The projection of the velocity vector field $\overrightarrow V $ on
the eigenbasis is written according to Corollary 2:

\[
\frac{d\delta \vec {X}}{dt} = \overrightarrow V _T \approx
\overrightarrow V = \sum\limits_{i = 1}^n {a_i } \lambda _i
\overrightarrow {Y_{\lambda _i } } = \alpha \overrightarrow
{Y_{\lambda _1 } } + \beta \overrightarrow {Y_{\lambda _2 } }
\]

where $\alpha $ and $\beta $ represent coefficients depending
explicitly on the co-ordinates of space and implicitly on time and
where $\overrightarrow {Y_{\lambda _1 } } $ represents the ``fast''
eigenvector and $\overrightarrow {Y_{\lambda _2 } } $ the ``slow''
eigenvector. The existence of an evanescent mode in the vicinity of
the \textit{slow manifold} implies according to Tihonv's theorem
[1952]: $\alpha \ll 1$ We deduce:\\

\textbf{Proposition A-1:} \textit{A necessary and sufficient
condition of obtaining the slow manifold equation of a two
dimensional dynamical system is that its velocity vector field
}$\overrightarrow V $\textit{ is collinear to the slow eigenvector
}$\overrightarrow {Y_{\lambda _2 } } $\textit{ associated in the
slow eigenvalue }$\lambda _2 $\textit{ of the functional jacobian of
the tangent linear system. That is to say:}

\[
\overrightarrow V \approx \beta \overrightarrow {Y_{\lambda _2 } }
\qquad (A-10)
\]

While using this collinearity condition, the equation constituting
the first order approximation in $\varepsilon $ of the \textit{slow
manifold} of a two dimensional dynamical system is written:

\[
\overrightarrow V \wedge \overrightarrow {Y_{\lambda _2 } } = \vec
{0}\mbox{ } \Leftrightarrow \mbox{ }\left( {\frac{\partial
g}{\partial x}} \right)\left( {\frac{dx}{dt}} \right) - \left(
{\lambda _2 - \frac{\partial g}{\partial y}} \right)\left(
{\frac{dy}{dt}} \right) = 0 \qquad (A-11)
\]

\subsubsection*{Slow manifold equation of a three
dimensional dynamical system}

Let's consider a dynamical system defined under the same conditions
as (\ref{eq1}) or (\ref{eq2}). The eigenvectors associated in the
eigenvalues of the functional jacobian of the \textit{tangent linear
system} are written:

\[
\overrightarrow {Y_{\lambda _i } } \left( {{\begin{array}{*{20}c}
 {\frac{1}{\varepsilon }\frac{\partial f}{\partial y}\frac{\partial
g}{\partial z} + \frac{1}{\varepsilon }\frac{\partial f}{\partial
z}\left(
{\lambda _i - \frac{\partial g}{\partial y}} \right)} \hfill \\
 {\frac{1}{\varepsilon }\frac{\partial f}{\partial z}\frac{\partial
g}{\partial x} + \frac{\partial g}{\partial z}\left( {\lambda _i -
\frac{1}{\varepsilon }\frac{\partial f}{\partial x}} \right)} \hfill \\
 { - \frac{1}{\varepsilon }\frac{\partial f}{\partial y}\frac{\partial
g}{\partial x} + \left( {\lambda _i - \frac{1}{\varepsilon
}\frac{\partial f}{\partial x}} \right)\left( {\lambda _i -
\frac{\partial g}{\partial y}}
\right)} \hfill \\
\end{array} }} \right) \qquad(A-12)
\]

with

\begin{center}
i =1, 2, 3
\end{center}

The projection of the velocity vector field $\overrightarrow V $ on
the eigenbasis is written according to Corollary 2:

\[
\frac{d\delta \vec {X}}{dt} = \overrightarrow V _T \approx
\overrightarrow V = \sum\limits_{i = 1}^n {a_i } \lambda _i
\overrightarrow {Y_{\lambda _i } } = \alpha \overrightarrow
{Y_{\lambda _1 } } + \beta \overrightarrow {Y_{\lambda _2 } } +
\delta \overrightarrow {Y_{\lambda _3 } }
\]

where $\alpha$, $\beta$ and $\delta $ represent coefficients
depending explicitly on the co-ordinates of space and implicitly on
time and where $\overrightarrow {Y_{\lambda _1 } } $ represents the
``fast'' eigenvector and $\overrightarrow {Y_{\lambda _2 } } $,
$\overrightarrow {Y_{\lambda _3 } } $ the ``slow'' eigenvectors. The
existence of an evanescent mode in the vicinity of the \textit{slow
manifold} implies according to Tihonv's theorem
[1952]: $\alpha \ll 1$ We deduce:\\

\textbf{Proposition~A-2:} \textit{A necessary and sufficient
condition of obtaining the slow manifold equation of a three
dimensional dynamical system is that its velocity vector field
}$\overrightarrow V $\textit{ is coplanar to the slow eigenvectors
}$\overrightarrow {Y_{\lambda _2 } } $\textit{ and }$\overrightarrow
{Y_{\lambda _3 } } $\textit{ associated in the slow eigenvalue
}$\lambda _2 $\textit{ and }$\lambda _3 $\textit{ of the functional
jacobian of the tangent linear system. That is to say:}

\[
\overrightarrow V \approx \beta \overrightarrow {Y_{\lambda _2 } } +
\delta \overrightarrow {Y_{\lambda _3 } } \qquad (A-13)
\]

While using this coplanarity condition, the equation constituting
the first order approximation in $\varepsilon $ of the \textit{slow
manifold} of a three dimensional dynamical system is written:

\[
\det \left( {\overrightarrow V ,\overrightarrow {Y_{\lambda _2 } }
,\overrightarrow {Y_{\lambda _3 } } } \right) = 0\mbox{ }
\Leftrightarrow \mbox{ }\overrightarrow V .\left( {\overrightarrow
{Y_{\lambda _2 } } \wedge \overrightarrow {Y_{\lambda _3 } } }
\right) = 0 \qquad(A-14)
\]

\subsection*{Orthogonality condition}

\subsubsection*{Slow manifold equation of a two dimensional
dynamical system}

Let's consider a dynamical system defined under the same conditions
as (\ref{eq1}) or (\ref{eq2}). The eigenvectors associated in the
eigenvalues of the transposed functional jacobian of the
\textit{tangent linear system} are written:

\[
{ }^t\overrightarrow {Y_{\lambda _i } } \left(
{{\begin{array}{*{20}c}
 {\lambda _i - \frac{\partial g}{\partial y}} \hfill \\
 {\frac{1}{\varepsilon }\frac{\partial f}{\partial y}} \hfill \\
\end{array} }} \right) \qquad(A-15)
\]

with

\begin{center}
i =1, 2
\end{center}

${ }^t\overrightarrow {Y_{\lambda _1 } } $represents the ``fast''
eigenvector on the left associated in the dominant eigenvalue, i.e.,
the largest eigenvalue in absolute value and ${ }^t\overrightarrow
{Y_{\lambda _2 } } $ is the ``slow'' eigenvector
on the left.\\

But since according to Rossetto \textit{et al.} [1998], the velocity
vector field $\overrightarrow V $ is perpendicular to the ``fast''
eigenvector on the left ${ }^t\overrightarrow {Y_{\lambda _1 } } $,
we deduce:\\

\textbf{Proposition A-3:} \textit{A necessary and sufficient
condition of obtaining the slow manifold equation of a two
dimensional dynamical system is that its velocity vector field
}$\overrightarrow V $\textit{ is perpendicular to the fast
eigenvector }${ }^t\overrightarrow {Y_{\lambda _1 } } $\textit{on
the left associated in the fast eigenvalue }$\lambda _1 $\textit{ of
the transposed functional jacobian of the tangent linear system.
That is to say:}
\[
\overrightarrow V \bot { }^t\overrightarrow {Y_{\lambda _1 } }
\qquad (A-16)
\]

While using this orthogonality condition, the equation constituting
the first order approximation in $\varepsilon $ of the \textit{slow
manifold} of a two dimensional dynamical system is written:

\[
\overrightarrow V \cdot { }^t\overrightarrow {Y_{\lambda _1 } } =
0\mbox{ } \Leftrightarrow \mbox{ }\left( {\lambda _1 -
\frac{\partial g}{\partial y}} \right)\left( {\frac{dx}{dt}} \right)
+ \left( {\frac{1}{\varepsilon }\frac{\partial f}{\partial y}}
\right)\left( {\frac{dy}{dt}} \right) = 0 \qquad(A-17)
\]

\newpage

\subsubsection*{Slow manifold equation of a three
dimensional dynamical system}

Let's consider a dynamical system defined under the same conditions
as (\ref{eq1}) or (\ref{eq2}). The eigenvectors associated in the
eigenvalues of the transposed functional jacobian of the
\textit{tangent linear system} are written:

\[
{ }^t\overrightarrow {Y_{\lambda _i } } \left(
{{\begin{array}{*{20}c}
 {\frac{\partial g}{\partial x}\frac{\partial h}{\partial y} +
\frac{\partial h}{\partial x}\left( {\lambda _i - \frac{\partial
g}{\partial
y}} \right)} \hfill \\
 {\frac{1}{\varepsilon }\frac{\partial f}{\partial y}\frac{\partial
h}{\partial x} + \frac{\partial h}{\partial y}\left( {\lambda _i -
\frac{1}{\varepsilon }\frac{\partial f}{\partial x}} \right)} \hfill \\
 { - \frac{1}{\varepsilon }\frac{\partial f}{\partial y}\frac{\partial
g}{\partial x} + \left( {\lambda _i - \frac{1}{\varepsilon
}\frac{\partial f}{\partial x}} \right)\left( {\lambda _i -
\frac{\partial g}{\partial y}}
\right)} \hfill \\
\end{array} }} \right) \qquad(A-18)
\]

with

\begin{center}
i =1, 2, 3
\end{center}

${ }^t\overrightarrow {Y_{\lambda _1 } } $represents the ``fast''
eigenvector on the left associated in the dominant eigenvalue, i.e.,
the largest eigenvalue in absolute value and ${ }^t\overrightarrow
{Y_{\lambda _2 } } $, ${ }^t\overrightarrow {Y_{\lambda _3 } } $ are
the ``slow'' eigenvectors on the left.
\\
But since according to Rossetto \textit{et al.} [1998], the velocity
vector field $\overrightarrow V $ is perpendicular to the ``fast''
eigenvector on the left ${ }^t\overrightarrow {Y_{\lambda _1 } } $,
we deduce:\\

\textbf{Proposition~A-4:} \textit{A necessary and sufficient
condition of obtaining the slow manifold equation of a three
dimensional dynamical system is that its velocity vector field
}$\overrightarrow V $\textit{ is perpendicular to the fast
eigenvector }${ }^t\overrightarrow {Y_{\lambda _1 } } $\textit{on
the left associated in the fast eigenvalue }$\lambda _1 $\textit{ of
the transposed functional jacobian of the tangent linear system.
That is to say:}

\[
\overrightarrow V \bot { }^t\overrightarrow {Y_{\lambda _1 } }
\qquad(A-19)
\]

While using this orthogonality condition, the equation constituting
the first order approximation in $\varepsilon $ of the \textit{slow
manifold} of a three dimensional dynamical system is written:

\[
\overrightarrow V \bot { }^t\overrightarrow {Y_{\lambda _1 } }
\mbox{ } \Leftrightarrow \mbox{ }\overrightarrow V \cdot {
}^t\overrightarrow {Y_{\lambda _1 } } = 0 \qquad(A-20)
\]

\newpage

\subsection*{Equivalence of both conditions}

\textbf{Proposition~A-5:} \textit{Both necessary and sufficient
collinearity / coplanarity and orthogonality conditions providing
the slow manifold equation are equivalent.}\\

\subsubsection*{Proof of the Proposition 5 in dimension two}

In dimension two, the \textit{slow manifold} equation may be
obtained while considering that the velocity vector field
$\overrightarrow V $ is:\\

- either collinear to the ``slow'' eigenvector $\overrightarrow
{Y_{\lambda _2 } } $

- either perpendicular to the ``fast'' eigenvector on the left ${
}^t\overrightarrow {Y_{\lambda _1 } } $\\

There is equivalence between both conditions provided that the
``fast'' eigenvector on the left ${ }^t\overrightarrow {Y_{\lambda
_1 } } $ is orthogonal to the ``slow'' eigenvector $\overrightarrow
{Y_{\lambda _2 } } $. Both co-ordinates of these eigenvectors
defined in the above section make it possible to express their
scalar product:

\[
{ }^t\overrightarrow {Y_{\lambda _1 } } \cdot \overrightarrow
{Y_{\lambda _2 } } = \lambda _1 \lambda _2 - \frac{\partial
g}{\partial y}\left( {\lambda _1 + \lambda _2 } \right) + \left(
{\frac{\partial g}{\partial y}} \right)^2 + \left(
{\frac{1}{\varepsilon }\frac{\partial f}{\partial y}} \right)\left(
{\frac{\partial g}{\partial x}} \right)
\]

While using the trace and determinant of the functional jacobian of
the \textit{tangent linear system}, we have:

\[
{ }^t\overrightarrow {Y_{\lambda _1 } } \cdot \overrightarrow
{Y_{\lambda _2 } } = \frac{1}{\varepsilon }\left( {\frac{\partial
f}{\partial x}\frac{\partial g}{\partial y} - \frac{\partial
f}{\partial y}\frac{\partial g}{\partial x}} \right) - \left(
{\frac{\partial g}{\partial y}} \right)\left( {\frac{1}{\varepsilon
}\frac{\partial f}{\partial x} + \frac{\partial g}{\partial y}}
\right) + \left( {\frac{\partial g}{\partial y}} \right)^2 + \left(
{\frac{1}{\varepsilon }\frac{\partial f}{\partial y}} \right)\left(
{\frac{\partial g}{\partial x}} \right) = 0
\]

So,

\[
{ }^t\overrightarrow {Y_{\lambda _1 } } \bot \overrightarrow
{Y_{\lambda _2 } } \qquad(A-21)
\]

Thus, collinearity and orthogonality conditions are completely
equivalent.

\newpage

\subsubsection*{Proof of the Proposition 5 in dimension three}

In dimension three, the \textit{slow manifold} equation may be
obtained while considering that the velocity vector field
$\overrightarrow V $ is:\\

- either coplanar to the ``slow'' eigenvectors $\overrightarrow
{Y_{\lambda _2 } } $ and $\overrightarrow {Y_{\lambda _3 } } $

- either perpendicular to the ``fast'' eigenvector on the left ${
}^t\overrightarrow {Y_{\lambda _1 } } $\\

There is equivalence between both conditions provided that the
``fast'' eigenvector on the left ${ }^t\overrightarrow {Y_{\lambda
_1 } } $ is orthogonal to the plane containing the ``slow''
eigenvectors $\overrightarrow {Y_{\lambda _2 } } $ and
$\overrightarrow {Y_{\lambda _3 } } $. While using the co-ordinates
of these eigenvectors defined in the above section, the sum and the
product (also the square of the sum and the product) of the
eigenvalues of the functional jacobian of the \textit{tangent linear
system}, the following equality is demonstrated:

\[
\overrightarrow {Y_{\lambda _2 } } \wedge \overrightarrow
{Y_{\lambda _3 } } = { }^t\overrightarrow {Y_{\lambda _1}} \qquad
(A-22)
\]

Thus, coplanarity and orthogonality conditions are completely
equivalent.

\begin{flushleft}
\textbf{Note:}
\end{flushleft}

In dimension three numerical studies shown significant differences
in the plot of the \textit{slow manifold} according to whether one
or the other of the both conditions (coplanarity (A-14) or
orthogonality (A-20)) were used. These differences come from the
fact that each one of these two conditions uses one or two
eigenvectors whose co-ordinates are expressed according to the
eigenvalues of the functional jacobian of the \textit{tangent linear
system}. Since these eigenvalues can be complex or real according to
their localization in the phase space the plot of the analytical
equation of the \textit{slow manifold} can be difficult even
impossible. Also to solve this problem it is necessary to make the
analytical equation of the \textit{slow manifold} independent of the
eigenvalues. This can be carried out by multiplying each equation of
the \textit{slow manifold} by one or two ``conjugated'' equations.
The equation obtained will be presented in each case (dimension two
and three) in the next section.

\newpage

\subsection*{Slow manifold equation independent of the eigenvectors}

\textbf{Proposition A-6:} \textit{Slow manifold equations of a
dynamical system obtained by the collinearity / coplanarity and
orthogonality conditions are equivalent. }

\subsubsection*{Proof of the Proposition 6 in dimension two}

In order to demonstrate the equivalence between the \textit{slow
manifold} equations obtained by each condition, they should be
expressed independently of the eigenvalues. So let's multiply each
equation (A-11) then (A-17) by its ``conjugated'' equation, i.e., an
equation in which the eigenvalue $\lambda _1 $ (resp. $\lambda _2 )$
is replaced by the eigenvalue $\lambda _2 $ (resp. $\lambda _1 )$.
Let's notice that the ``conjugated'' equation of the equation (A-11)
corresponds to the collinearity condition between the velocity
vector field $\overrightarrow V $ and the eigenvector
$\overrightarrow {Y_{\lambda _1 } } $. The product of the equation
(A-11) by its ``conjugated'' equation is written:

\[
\left( {\overrightarrow V \wedge \overrightarrow {Y_{\lambda _1 } }
} \right) \cdot \left( {\overrightarrow V \wedge \overrightarrow
{Y_{\lambda _2 } } } \right) = 0
\]

So, while using the trace and the determinant of the functional
jacobian of the \textit{tangent linear system} we have:

\[
\left( {\frac{\partial g}{\partial x}} \right)\left[ {\left(
{\frac{\partial g}{\partial x}} \right)\left( {\frac{dx}{dt}}
\right)^2 - \left( {\frac{1}{\varepsilon }\frac{\partial f}{\partial
x} - \frac{\partial g}{\partial y}} \right)\left( {\frac{dx}{dt}}
\right)\left( {\frac{dy}{dt}} \right) - \left( {\frac{1}{\varepsilon
}\frac{\partial f}{\partial y}} \right)\left( {\frac{dy}{dt}}
\right)^2} \right] = 0 \qquad (A-23)
\]

In the same manner, the product of the equation (A-17) by its
``conjugated'' equation which corresponds to the orthogonality
condition between velocity vector field $\overrightarrow V $ and the
eigenvector ${ }^t\overrightarrow {Y_{\lambda _2 } } $ is written:

\[
\left( {\overrightarrow V \cdot { }^t\overrightarrow {Y_{\lambda _1
} } } \right)\left( {\overrightarrow V \cdot { }^t\overrightarrow
{Y_{\lambda _2 } } } \right) = 0
\]

So, while using the trace and the determinant of the functional
jacobian of the \textit{tangent linear system} we have:

\[
\left( {\frac{1}{\varepsilon }\frac{\partial f}{\partial y}}
\right)\left[ {\left( {\frac{\partial g}{\partial x}} \right)\left(
{\frac{dx}{dt}} \right)^2 - \left( {\frac{1}{\varepsilon
}\frac{\partial f}{\partial x} - \frac{\partial g}{\partial y}}
\right)\left( {\frac{dx}{dt}} \right)\left( {\frac{dy}{dt}} \right)
- \left( {\frac{1}{\varepsilon }\frac{\partial f}{\partial y}}
\right)\left( {\frac{dy}{dt}} \right)^2} \right] = 0 \qquad (A-24)
\]

Both equations (A-23) and (A-24) are equal provided that:

\[
\left( {\frac{\partial g}{\partial x}} \right) \ne 0 \mbox{ and }
\left( {\frac{1}{\varepsilon }\frac{\partial f}{\partial y}} \right)
\ne 0
\]

These two last conditions are, according to the definition of a
dynamical system, satisfied because if they were not both
differential equations which constitute the system would be
completely uncoupled and it would not be a system anymore. Thus, the
equations obtained by the collinearity and orthogonality conditions
are equivalent.

\[
\left( {\overrightarrow V \wedge \overrightarrow {Y_{\lambda _1 } }
} \right) \cdot \left( {\overrightarrow V \wedge \overrightarrow
{Y_{\lambda _2 } } } \right) = 0\mbox{ } \Leftrightarrow \mbox{
}\left( {\overrightarrow V \cdot { }^t\overrightarrow {Y_{\lambda _1
} } } \right)\left( {\overrightarrow V \cdot { }^t\overrightarrow
{Y_{\lambda _2 } } } \right) = 0 \qquad(A-25)
\]

Equations (A-23) and (A-24) provide the \textit{slow manifold}
equation of a two dimensional dynamical system independently of the
eigenvalues of the functional jacobian of the \textit{tangent linear
system}. In order to express them, we adopt the following notations
for:\\

- the velocity vector field

\[
\overrightarrow V \left( {{\begin{array}{*{20}c}
 \dot {x} \hfill \\
 \dot {y} \hfill \\
\end{array} }} \right)
\]

- the functional jacobian

\[
J = \left( {{\begin{array}{*{20}c}
 a \hfill & b \hfill \\
 c \hfill & d \hfill \\
\end{array} }} \right)
\]

- the eigenvectors co-ordinates (A-9)

\[
\overrightarrow {Y_{\lambda _i } } \left( {{\begin{array}{*{20}c}
 {\lambda _i - d} \hfill \\
 c \hfill \\
\end{array} }} \right)
\]

\begin{center}
with i = 1, 2
\end{center}

Equations (A-23) and (A-24) providing the \textit{slow manifold}
equation of a two dimensional dynamical system independently of the
eigenvalues of the functional jacobian of the \textit{tangent linear
system} are then written:

\[
c\dot {x}^2 - \left( {a - d} \right)\dot {x}\dot {y} - b\dot {y}^2 =
0 \qquad(A-26)
\]

\[
\phi = \sum\limits_{i,j = 0}^2 {\alpha _{ij} \dot {x}^i\dot {y}^j} =
0 \mbox{ with } \alpha _{ij} = \left\{ {{\begin{array}{*{20}c}
 { = 0\mbox{ if }i + j \ne 2} \hfill \\
 { \ne 0\mbox{ if }i + j = 2} \hfill \\
\end{array} }} \right.\qquad (A-27)
\]

\begin{table}[htbp]
\begin{center}
\begin{tabular}{|c|}
\hline
\textbf{Slow manifold analytical equation of a two dimensional dynamical system}   \\
\hline $\phi = \sum\limits_{i,j = 0}^2 {\alpha _{ij} \dot {x}^i\dot
{y}^j} = 0 $ with $ \alpha _{ij} = \left\{ {{\begin{array}{*{20}c}
 { = 0\mbox{ if }i + j \ne 2} \hfill \\
 { \ne 0\mbox{ if }i + j = 2} \hfill \\
\end{array} }} \right.$ \\
\hline
$\alpha _{20} = c$ \\
\hline
$\alpha _{11} = - \left( {a - d} \right)$ \\
\hline
$\alpha _{02} = - b$ \\
\hline
\end{tabular}
\end{center}
\caption{Slow manifold analytical equation of a two dimensional
dynamical system} \label{tab2}
\end{table}

\subsubsection*{Proof of the Proposition 6 in dimension three}

In order to demonstrate the equivalence between the \textit{slow
manifold} equations obtained by each condition, the same step as
that exposed in the above section is applied. The \textit{slow
manifold} equation should be expressed independently of the
eigenvalues. In dimension three, each equation (A-14) then (A-20)
must be multiplied by two ``conjugated'' equations obtained by
circular shifts of the eigenvalues. Let's notice that the first of
the ``conjugated'' equations of the equation (A-14) corresponds to
the coplanarity condition between the velocity vector field
$\overrightarrow V $ and the eigenvectors $\overrightarrow
{Y_{\lambda _1 } } $and $\overrightarrow {Y_{\lambda _2 } } $, the
second corresponds to the coplanarity condition between the velocity
vector field $\overrightarrow V $ and the eigenvectors
$\overrightarrow {Y_{\lambda _1 } } $and $\overrightarrow
{Y_{\lambda _3 } } $. The product of the equation (A-14) by its
``conjugated'' equation is written:
\[
\left[ {\overrightarrow V \cdot \left( {\overrightarrow {Y_{\lambda
_1 } } \wedge \overrightarrow {Y_{\lambda _2 } } } \right)} \right]
\cdot \left[ {\overrightarrow V \cdot \left( {\overrightarrow
{Y_{\lambda _2 } } \wedge \overrightarrow {Y_{\lambda _3 } } }
\right)} \right] \cdot \left[ {\overrightarrow V \cdot \left(
{\overrightarrow {Y_{\lambda _1 } } \wedge \overrightarrow
{Y_{\lambda _3 } } } \right)} \right] = 0 \qquad (A-28)
\]

In the same manner, the product of the equation (A-20) by its
``conjugated'' equation which corresponds to the orthogonality
condition between the velocity vector field $\overrightarrow V $ and
the eigenvector ${ }^t\overrightarrow {Y_{\lambda _2 } } $ and, the
orthogonality condition between the velocity vector field
$\overrightarrow V $ and the eigenvector ${ }^t\overrightarrow
{Y_{\lambda _3 } } $ is written:

\[
\left( {\overrightarrow V \cdot { }^t\overrightarrow {Y_{\lambda _1
} } } \right)\left( {\overrightarrow V \cdot { }^t\overrightarrow
{Y_{\lambda _2 } } } \right)\left( {\overrightarrow V \cdot {
}^t\overrightarrow {Y_{\lambda _3 } } } \right) = 0 \qquad (A-29)
\]

By using Eq. (A-22) and all the circular shifts which result from
this we demonstrate that Eqs. (A-28) and (A-29) are equal. Thus, the
equations obtained by the coplanarity and orthogonality conditions
are equivalent.

\newpage

\[
\left[ {\overrightarrow V \cdot \left( {\overrightarrow {Y_{\lambda
_1 } } \wedge \overrightarrow {Y_{\lambda _2 } } } \right)} \right]
\cdot \left[ {\overrightarrow V \cdot \left( {\overrightarrow
{Y_{\lambda _2 } } \wedge \overrightarrow {Y_{\lambda _3 } } }
\right)} \right] \cdot \left[ {\overrightarrow V \cdot \left(
{\overrightarrow {Y_{\lambda _1 } } \wedge \overrightarrow
{Y_{\lambda _3 } } } \right)} \right] = 0
\]

\[
\Leftrightarrow \qquad(A-30)
\]

\[
\left( {\overrightarrow V \cdot { }^t\overrightarrow {Y_{\lambda _1
} } } \right)\left( {\overrightarrow V \cdot { }^t\overrightarrow
{Y_{\lambda _2 } } } \right)\left( {\overrightarrow V \cdot {
}^t\overrightarrow {Y_{\lambda _3 } } } \right) = 0
\]

Equations (A-28) and (A-29) provide the \textit{slow manifold}
equation of a three dimensional dynamical system independently of
the eigenvalues of the functional jacobian of the \textit{tangent
linear system}. In order to express them, we adopt the following
notations for:\\

- the velocity vector field

\[
\overrightarrow V \left( {{\begin{array}{*{20}c}
 \dot {x} \hfill \\
 \dot {y} \hfill \\
 \dot {z} \hfill \\
\end{array} }} \right)
\]

- the functional jacobian

\[
J = \left( {{\begin{array}{*{20}c}
 a \hfill & b \hfill & c \hfill \\
 d \hfill & e \hfill & f \hfill \\
 g \hfill & h \hfill & i \hfill \\
\end{array} }} \right)
\]

- the ``slow'' eigenvectors co-ordinates (A-12)

\[
\overrightarrow {Y_{\lambda _i } } \left( {{\begin{array}{*{20}c}
 {bf + c\left( {\lambda _i - e} \right)} \hfill \\
 {cd + f\left( {\lambda _i - a} \right)} \hfill \\
 { - bd + \left( {\lambda _i - a} \right)\left( {\lambda _i - e} \right)}
\hfill \\
\end{array} }} \right)
\]

\begin{center}
with i = 2, 3
\end{center}

- the ``fast'' eigenvectors co-ordinates on the left (A-18)

\[
{ }^t\overrightarrow {Y_{\lambda _1 } } \left(
{{\begin{array}{*{20}c}
 {hd + g\left( {\lambda _1 - e} \right)} \hfill \\
 {bg + h\left( {\lambda _1 - a} \right)} \hfill \\
 { - bd + \left( {\lambda _1 - a} \right)\left( {\lambda _1 - e} \right)}
\hfill \\
\end{array} }} \right)
\]

Starting from the coplanarity condition (A-14) and while replacing
the eigenvectors by their co-ordinates (A-12) and while removing all
the eigenvalues $\lambda _2 $ and $\lambda _3 $ thanks to the sum
and the products (and also the square of the sum and the products)
of the eigenvalues of the functional jacobian of the \textit{tangent
linear} \textit{system}, we obtain the following equation:

\[
A_1 \dot {x} - B_1 \dot {y} + C\dot {z} = 0 \qquad (A-31)
\]

with

\[
A_1 = f\lambda _1^2 - \left( {ef + if + cd} \right)\lambda _1 + efi
+ cdi - cfg - f^2h
\]

\[
B_1 = c\lambda _1^2 - \left( {ac + ic + bf} \right)\lambda _1 + aci
+ bfi - c^2g - cfh \qquad (A-32)
\]

\[
C = bf^2 - c^2d + cf\left( {a - e} \right)
\]

The equation (A-31) is absolutely identical to that which one would
have obtained by the orthogonality condition (A-20). Let's multiply
the equation (A-31) by its ``conjugated'' equations in $\lambda _2 $
and $\lambda _3 $, i.e., by $\left( {A_2 \dot {x} - B_2 \dot {y} +
C\dot {z}} \right)$ and by $\left( {A_3 \dot {x} - B_3 \dot {y} +
C\dot {z}} \right)$. The coefficients $A_i $, $B_i $ are obtained by
replacing in the equations (A-32) the eigenvalue $\lambda _1 $ by
the eigenvalue $\lambda _2 $ then by the eigenvalue $\lambda _3 $
respectively for i = 2, 3. We obtain:

\[
\left( {A_1 \dot {x} - B_1 \dot {y} + C\dot {z}} \right)\left( {A_2
\dot {x} - B_2 \dot {y} + C\dot {z}} \right)\left( {A_3 \dot {x} -
B_3 \dot {y} + C\dot {z}} \right) = 0 \qquad (A-33)
\]

So,

\[
\phi = \sum\limits_{i,j,k = 0}^3 {\alpha _{ijk} \dot {x}^i\dot
{y}^j\dot {z}^k} = 0 \mbox{ with } \alpha _{ijk} = \left\{
{{\begin{array}{*{20}c}
 { = 0\mbox{ if }i + j + k \ne 3} \hfill \\
 { \ne 0\mbox{ if }i + j + k = 3} \hfill \\
\end{array} }} \right.\qquad (A-34)
\]

By developing this expression we obtain a polynomial comprising
terms comprising the sum and product of the eigenvalues and also of
the square of the sum and the product of the eigenvalues and which
are directly connected to the elements of the functional jacobian
matrix of the \textit{tangent linear system}. The equation obtained
is the result of a demonstration (available by request to the
authors) which establishes a relation between the coefficients of
this polynomial and the elements of the functional jacobian matrix
of the \textit{tangent linear system}.

\begin{table}[htbp]
\begin{center}
\begin{tabular}
{|c|} \hline
\textbf{Slow manifold analytical equation of a three dimensional dynamical system} \par   \\
\hline $\phi = \sum\limits_{i,j,k = 0}^3 {\alpha _{ijk} \dot
{x}^i\dot {y}^j\dot {z}^k} = 0 $ with $ \alpha _{ijk} = \left\{
{{\begin{array}{*{20}c}
 { = 0\mbox{ if }i + j + k \ne 3} \hfill \\
 { \ne 0\mbox{ if }i + j + k = 3} \hfill \\
\end{array} }} \right.$ \\
\hline
$\alpha _{300} = d^2h + dgi - fg^2 - dge$ \\
\hline
$\alpha _{030} = ch^2 + abh - b^2g - ibh$ \\
\hline
$\alpha _{003} = c^2d + cfe - bf^2 - cfa$ \\
\hline
$\alpha _{210} = bdg + aeg - e^2g + cg^2 - 2adh - 2fgh + deh - agi + egi + dhi$ \\
\hline
$\alpha _{120} = - abg + 2beg + a^2h - fh^2 - bdh - aeh + 2cgh - bgi - ahi + ehi$ \\
\hline
$\alpha _{201} = - bd^2 + ade - cdg + 2afg + 2dfh - efg - adi - dei - fgi + di^2$ \\
\hline
$\alpha _{102} = acd + cde - a^2f - 2bdf + aef + cfg + f^2h - 2cdi + afi - efi$ \\
\hline
$\alpha _{021} = b^2d - abe - 2bcg - 2ceh + ach + bfh + abi + bei + chi - bi^2$ \\
\hline
$\alpha _{012} = 2bcd - ace + ce^2 - abf - bef - c^2g - cfh + aci - cei + 2bfi$ \\
\hline
$\alpha _{111} = abd - a^2e - bde + ae^2 - acg + 3bfg - 3cdh + efh + a^2i - e^2i + cgi - fhi - ai^2 + ei^2$ \\
\hline
\end{tabular}
\end{center}
\caption{Slow manifold analytical equation of a three dimensional
dynamical system} \label{tab3}
\end{table}

The expression (A-34) represents the \textit{slow manifold }equation
of a three dimensional dynamical system independently of the
eigenvalues of the functional jacobian of the \textit{tangent linear
system}. Both expression (A-27) and (A-34) are also available at the
address: \href{http://ginoux.univ-tln.fr}{http://ginoux.univ-tln.fr}

\end{document}